%% file: wienerbasis-part1.tex
\documentclass{article}
\usepackage{amsmath,amssymb,graphicx,bbm,theorem}
\pdfoutput=1
\usepackage{tikz,verbatim}
\usepackage{pgfmath}
\usetikzlibrary{decorations}
\usepgflibrary{decorations.pathreplacing}
\usepgflibrary{shapes.misc}
\usetikzlibrary{arrows}

\makeatletter

\newcommand{\Rmnum}[1]{\expandafter\@slowromancap\romannumeral #1@}
\makeatother

\usepackage[footnotesize,bf]{caption}  
\usepackage{array} 

\usepackage{amsmath,amssymb,graphicx,bbm,theorem}
\usepackage{verbatim}
\usepackage{epsfig}
\usepackage{wrapfig}
\usepackage{latexsym,amsfonts,amscd,subfigure}
\usepackage{changebar}
\usepackage{enumerate}

\newcommand{\assign}{:=}
\newcommand{\asterisk}{*}

\newcommand{\mathd}{\mathrm{d}}
\newcommand{\nin}{\not\in}
\newcommand{\nonesep}{}
\newcommand{\tmem}[1]{{\em #1\/}}
\newcommand{\tmop}[1]{\ensuremath{\operatorname{#1}}}
\newcommand{\tmscript}[1]{\text{\scriptsize{$#1$}}}
\newcommand{\tmtextit}[1]{{\itshape{#1}}}
\newenvironment{proof}{\noindent\textbf{Proof\ }}{\hspace*{\fill}$\Box$\medskip}
\newtheorem{corollary}{Corollary}
\newtheorem{lemma}{Lemma}
\newtheorem{proposition}{Proposition}
{\theorembodyfont{\rmfamily}\newtheorem{remark}{Remark}}
\newtheorem{theorem}{Theorem}

\begin{document}

\title{A Generalization of the Wiener Rational Basis \\
Functions on Infinite Intervals\\
Part I -- Derivation and Properties}
\author{Akil C. Narayan\thanks{Division of Applied Mathematics; Brown University; 182 George
Street, Box F; Providence, RI 02912; {\tt anaray@dam.brown.edu}} 
\and Jan S.
Hesthaven\thanks{Division of Applied Mathematics; Brown University; {\tt
Jan.Hesthaven@brown.edu}}}
\date{}

\maketitle

\begin{abstract}
  \input{part1Content/abstract}
\end{abstract}

\section{Introduction}
\input{part1Content/introduction}

\section{Derivation of the basis}
\input{part1Content/derivation}

\section{Fourier-Derived Basis Properties}
\input{part1Content/fourierprop}

\section{Jacobi-Derived Basis Properties}
\input{part1Content/jacobiprop}

\section{The Semi-Infinite Interval}
\input{part1Content/semiinfinite.tex}

\section{Alternative Methods}
\input{part1Content/mjpoly.tex}

\section{Conclusion}
\input{part1Content/conclusion}

{\newpage}
\appendix
\input{part1Content/appendix}

\bibliographystyle{amsplain}
\bibliography{wienerbasis}

\end{document}

%% file: part1Content/abstract.tex
  \noindent We formulate and derive a generalization of an orthogonal
  rational-function basis for spectral expansions over the infinite or
  semi-infinite interval. The original functions, first presented by Wiener
  {\cite{wiener1949}}, are a mapping and weighting of the Fourier basis to the
  infinite interval. By identifying the Fourier series as a biorthogonal
  composition of Jacobi polynomials/functions, we are able to define
  generalized Fourier series' which, appropriately mapped to the whole real line
  and weighted, form a generalization of Wiener's basis functions. It is known
  that the original Wiener rational functions inherit sparse Galerkin matrices
  for differentiation, and can utilize the fast Fourier transform (FFT) for
  computation of the modal coefficients.  We show that the generalized basis
  sets also have a sparse differentiation matrix and we discuss connection
  problems, which are necessary theoretical developments for application of the
  FFT.

%% file: part1Content/introduction.tex
\label{sec:intro}The approximation of a function by a finite sum of 
basis functions has long been a hallmark tool in numerical analysis. Over the
finite interval much is known about expansion properties and periodic Fourier
expansions or polynomial expansions are well-studied. On infinite intervals
there are complications due to the unbounded domain on which approximation is
necessary. Nevertheless many basis sets have been successfully investigated in
this case; the Hermite functions provide a suitable method for approximation when
it can be assumed that the function decays exponentially; for functions that
do not decay exponentially, the so-called mapped Chebyshev rational functions
can fill the void and open up the possibility for utilizing the fast Fourier
transform (FFT); additionally, a Fourier basis mapped to the real line has been
explored and provides an additional method for function approximation over the
infinite interval. This last basis set serves as an inspiration for the family
of basis sets proposed in this paper.

Despite the available methods for function approximation over the infinite
interval, there are shortcomings. The Hermite
functions/polynomials do not admit an FFT exploitation and have problems
approximating functions that do not decay exponentially (which is to say, most
functions). However, the solutions to differential equations have been
relatively successful by Hermite approximation and in some cases give superior
approximations when compared to a Chebyshev (mapped or truncated)
approximation {\cite{boyd1984}}. The Whittaker cardinal interpolant functions
{\cite{whittaker1915}}, or Sinc functions, provide a remarkably simple method
to approximate a function with known equispaced evaluations. The drawback is a
relatively small class of functions for which such an expansion is complete.
However, the ease of applying Sinc methods has led to a great number of
applications {\cite{lund1992}}. The Chebyshev rational functions
{\cite{boyd1982}}, {\cite{boyd1990a}} are robust with respect to the
deficiencies of the Hermite and Sinc bases, but they have some disadvantages
compared with the generalized Wiener basis we will derive. 

On the semi-infinite interval Laguerre polynomial/function expansions are the
classical approximation technique {\cite{abramowitz1972}}, but these techniques
suffer from the same problems as Hermite expansions. An alternative technique
involves mapping Jacobi polynomials to the infinite interval {\cite{boyd1987b}}.
This mapping technique makes it possible to accurately approximate
algebraically-decaying functions on the semi-infinite interval, but introduces
some computational issues for the solution to differential equations. The
generalized Wiener basis can be employed on the semi-infinite interval; this
results in a basis set that is also a mapped Jacobi polynomial methods. However,
the Wiener mapping is very different from that presented in the literature, and
therefore acts as a competitor to these existing techniques. 

Our generalized basis is inspired by a collection of
orthogonal and complete functions originally proposed by Wiener
{\cite{wiener1949}}. He introduces the functions
\begin{equation}
  \label{eq:wiener-functions} \phi_n (x) = \frac{(1 - i \nonesep
  x)^n}{\sqrt{\pi} (1 + i \nonesep x)^{n + 1}}, \hspace{0.5cm} n \in
  \mathbbm{N}_0
\end{equation}
as Fourier transforms of the Laguerre functions. He furthermore shows that
these functions are orthogonal under the $L^2$ conjugate inner product.
Higgins {\cite{higgins1977}} expands this result by presenting the functions
$\phi_n$ along with their complex conjugates as a complete system in $L^2$.
Following this, others have followed up on these functions by applying them to
the solution of differential equations {\cite{christov1982}}, {\cite{cain1984}}. We note that the
functions $\phi_n (x)$ presented above have magnitude that decays like
$\frac{1}{x}$ as $|x| \rightarrow \infty$. We will generalize the above
functions so that they have decay $\frac{1}{x^s}$ for any $s > \frac{1}{2}$.
The ability to choose the rate of decay of the basis set is an advantage if
such information is present about the nature of the function to be
approximated or the differential equation to be solved (e.g. {\cite{klaus2006}},
{\cite{magyari2004}}). Furthermore, we will
show that this basis admits sparse Galerkin matrices and that the fast Fourier
Transform can be used in certain cases to evaluate and manipulate the series.

This paper is concerned with the derivation and theoretical properties of the
generalized Wiener rational function basis. Computational considerations,
numerical examples, and comparisons with existing basis sets are presented in
a second part. The outline of this paper is as follows.
In Section \ref{sec:derivation} we formulate and derive the basis, which is
heavily based upon a generalized Fourier series. Section \ref{sec:prop}
follows with some properties of the basis functions based on their close
relationship to the canonical Fourier basis, and Section \ref{sec:propjacobi}
concerns the properties that can be derived from the relation to Jacobi
polynomials. In Section \ref{sec:semiinfinite} we discuss how the Wiener basis
set may be used to approximate functions on the semi-infinite interval. Finally,
we briefly present mapped Jacobi polynomials as an alternative method in Section
\ref{sec:mjpoly} and summarize and present an outlook in Section
\ref{sec:conclusion} for Part \Rmnum{2}, dealing with numerical issues.

%% file: part1Content/derivation.tex
\label{sec:derivation}We begin by stating the major goals and the path we will
take in accomplishing those goals. We seek a collection of
$L^2$-orthogonal and complete basis functions whose domain is the entire
real line. In addition, we desire the ability to specify a parameter $s >
\frac{1}{2}$ that will denote the polynomial decay at $\pm \infty$ of each of
the the basis functions.

Drawing inspiration from Wiener and his orthogonal basis functions, we seek a
collection of functions $\phi_k^{(s)} (x)$ for $x \in \mathbbm{R}$ and $k \in
\mathbbm{Z}$ such that $\left\{ \phi_k^{(s)} \right\}_{k \in \mathbbm{Z}}$ is
a complete, orthogonal system for any valid $s$. Our method relies on the
observation that the functions (\ref{eq:wiener-functions}) are weighted
maps of the canonical Fourier basis $e^{i \nonesep n \nonesep \theta}$ for
$\theta \in [0, 2 \pi]$ (see e.g. {\cite{boyd1990a}}, {\cite{weidemann1992}}).
We will first generalize the Fourier basis on $[0, 2 \pi]$ so that it will
have the properties we desire on the infinite interval; we will then map the
generalized Fourier basis to the real line and weight it accordingly to
achieve the desired rate of decay.

\subsection{Notation and setup}

\label{sec:derivation-notation}We shall reserve the variables $x, z, \theta,$
and $r$ as independent variables on certain domains and list the domains
and transformations in Table \ref{tab:xforms}. The variable $r \in [- 1, 1]$
is the standard interval over which the Jacobi polynomials are defined. The
interval $\theta \in [0, \pi]$ is the image of the $r$ interval under the map
$\theta = \arccos r$. The variable $z \in \mathbbm{T}^+$ is the upper-half of
the unit circle in the complex plane, and $x \in [0, \infty]$ is the positive
half of the extended real line.

In much of what follows we will mix notation and write expressions both in
terms of e.g. $r$ and $\theta$. It should then be understood that $r = r
(\theta)$ and/or $\theta = \theta (r)$. Furthermore, we will extend the
domains of $\theta, z,$ and $x$ to be $[- \pi, \pi]$, $\mathbbm{T}$, and
$\mathbbm{R}$, respectively later in the paper.

\begin{table}[tbp]
  
\newcolumntype{q}{>{$}c<{$}}
  \centering
  \renewcommand{\arraystretch}{2.0}
  \setlength{\doublerulesep}{0pt}
  \begin{tabular}{|q !{\vrule width 1.0pt} q|q|q|q|}\hline
       & x & z & \theta & r  \\\hline\hline
       x & x \in [0, \infty] & z = - \frac{x - i}{x + i} & \theta = 2 \arctan
       (x) & r = \frac{1 - x^2}{1 + x^2} \\\hline
       z & x = i \frac{1 - z}{1 + z}  & z \in \mathbbm{T}^+ & \theta = \arg z
       & r = \frac{1}{2} \left( z + \bar{z} \right) \\\hline
       \theta & x = \tan \left( \frac{\theta}{2} \right) & z = e^{i \theta} &
       \theta \in [0, \pi] & r = \cos \theta  \\\hline
       r & x = \sqrt{\frac{1 - r}{1 + r}} & z = e^{i \arccos r} & \theta =
       \arccos r & r \in [- 1, 1] \\\hline
   \end{tabular}
  \caption{Isomorphic transforms between different domains.}
\label{tab:xforms}
\end{table}

We denote $L^2 \left( A, B ; w) = L^2_w \left( A, B) \right. \right.$ the
space of square integrable functions $f : A \rightarrow B$ under the weight
$w$. We endow $L^2_w \left( A, B \right)$ with the conjugate bilinear inner
product; the notation for this inner product is $\left\langle \cdot, \cdot,
\right\rangle_w$. The omission of $w$ indicates the unit weight measure. The
norm on this space will be denoted $\left\| \cdot \right\|_w^{}$.  The following
weight functions will be used extensively in this article:

\begin{align*}
     w_r^{(\alpha, \beta)} (r) & = (1 - r)^{\alpha} (1 + r)^{\beta}\\
     &  \\
     w_{\theta}^{(\gamma, \delta)} (\theta) & = w_r^{(\delta, \gamma)} (r
     (\theta)) = (1 + \cos \theta)^{\gamma} (1 - \cos \theta)^{\delta}\\
     &  \\
     w_x^{(s, t)} (x) & = w_{\theta}^{(s, t)} (\theta (x)) = \frac{2^{s +
     t}}{(1 + x^2)^s}  \left( \frac{x^{2 t}}{(1 + x^2)^t} \right) .
\end{align*}
In addition, we will make use of a phase-shifted square root of $w_x^{(s, t)}$
and $w_{\theta}^{(\gamma, \delta)}$, which we define as:
\begin{equation}
  \label{eq:wx-npbranch} \sqrt[\asterisk]{w_x^{(s, t)} (x)} \assign
  \sqrt{w_x^{(s, t)}} \exp \left[ \frac{i (s + t)}{2} \left( \pi - \theta (x)
  \right) \right] = \frac{2^{\left( \frac{s + t}{2} \right)} x^t}{(x - i)^{s +
  t}}
\end{equation}
\begin{equation}
  \label{eq:wt-npbranch} \sqrt[\asterisk]{w_{\theta}^{(\gamma, \delta)}
  (\theta)} = \sqrt[\asterisk]{w_x^{(\gamma, \delta)} (x (\theta))} =
  2^{\left( \frac{\gamma + \delta}{2} \right)} \sin^{\delta} \left(
  \frac{\theta}{2} \right) \cos^{\gamma} \left( \frac{\theta}{2} \right) \exp
  \left[ \frac{i (\gamma + \delta)}{2} \left( \pi - \theta \right) \right]
\end{equation}

\subsection{Jacobi polynomials}

The classical Jacobi polynomials $P_n^{(\alpha, \beta)}$ are a family of
orthogonal polynomials {\cite{szego1975}} that have been used extensively in
many applications due to their ability to approximate general classes of
functions. They are a class of polynomials that encompass the Chebyshev,
Legendre, and Gegenbauer/ultraspheric polynomials. These polynomials will form
the building blocks for our generalization.

The Jacobi differential equation is
\begin{equation}
  \label{eq:jacobi-ode} \begin{array}{ll}
    (1 - r^2) \rho'' + \left[ \beta - \alpha - (\alpha + \beta + 2) r \right]
    \rho' + n (n + \alpha + \beta + 1) \rho = 0, & r \in [- 1, 1],
  \end{array}
\end{equation}
and for $\alpha, \beta > - 1$, $n \in \mathbbm{N}_0$ the only polynomial
solution $\rho = P^{(\alpha, \beta)}_n (x)$ is a polynomial of degree $n$. The
restriction $\alpha,\beta>-1$ is necessary to ensure integrability of the weight
and thus existence of an $L^2$-constant function solution. The
family of polynomials $\left\{ P_n^{(\alpha, \beta)} (x) \right\}_{n =
0}^{\infty}$ is complete and orthogonal in $L^2 \left( [- 1, 1], \mathbbm{R} ;
w_r^{(\alpha, \beta)} \right)$. We denote $h_n^{(\alpha, \beta)} = \left\|
P_n^{(\alpha, \beta)} \right\|_{w^{(\alpha, \beta)}_r}^2$, and define the
normalized polynomials as
\[ \tilde{P}^{(\alpha, \beta)}_n (r) = \frac{P^{(\alpha, \beta)}_n
   (r)}{\sqrt{h_n^{(\alpha, \beta)}}} . \]
The orthonormal Jacobi polynomials $\tilde{P}^{(\alpha, \beta)}_n$ will be
integral in the derivation of the Wiener rational function basis on the real
line. In addition, we require a minor generalization of Jacobi polynomials: we
perform a change of the dependent variable in (\ref{eq:jacobi-ode}) to
obtain:

\begin{lemma}
  \label{lemma:jacobi-functions}{\em (Jacobi Functions)}
  \noindent The Jacobi functions defined as
  \[ P^{(\alpha, \beta, a, b)}_n (r) = (1 - r)^a (1 + r)^b P_n^{(\alpha,
     \beta)} (r) \]
  satisfy the following properties:
  \begin{enumerate}
    \item $\left\{ P^{(\alpha, \beta, a, b)}_n (r) \right\}_{n \in
    \mathbbm{N}_0}$ are orthogonal and complete in $L^2\left( [-1,1],
    \mathbbm{R}; w_r^{(\alpha-2a,\beta-2b)}\right)$.
    
    \item The $P^{(\alpha, \beta, a, b)}_n (r)$ are eigenfunctions $\rho_n
    (r)$ of the Sturm-Liouville problem
    \begin{equation}
      \label{eq:slp-jacobi-functions} - \frac{\mathd}{\mathd r} \left[ p (r)
      \rho' (r) \right] + q (r) \rho (r) - \lambda_n w (r) \rho (r) = 0,
    \end{equation}
    which is defined by the parameters
    \begin{equation}
      \left. \label{eq:slp-jacobi-general-coeffs} \begin{array}{lll}
        p (r) & = & (1 - r)^{\alpha + 1 - 2 a} (1 + r)^{\beta + 1 - 2 b} \\
        &  & \\
        q (r) & = & \left[ a (\alpha - a) (1 - r)^{- 2} + b (\beta - b) (1 +
        r)^{- 2} \right] (1 - r)^{\alpha + 1 - 2 a} (1 + r)^{\beta + 1 - 2
        b}\\
        &  & \\
        w (r) & = & (1 - r)^{\alpha - 2 a} (1 + r)^{\beta - 2 b}\\
        &  & \\
        \lambda_n & = & n (n + \alpha + \beta + 1) - 2 ab + a (\beta + 1) + b
        (\alpha + 1)
      \end{array} \right\}
    \end{equation}
  \end{enumerate}
\end{lemma}

\noindent The proof is mathematically simple but algebraically tedious and we omit it.
We shall actually only require the result of Lemma
\ref{lemma:jacobi-functions} for $a = b = \frac{1}{2}$. Many of the results in
this paper require the use of numerous recurrence relations involving Jacobi
polynomials; these relations are given in Appendix \ref{app:recurrence},
equations (\ref{eq:opoly-3term})-(\ref{eq:jdiff}).

The idea behind the formation of the Jacobi functions introduced in Lemma
\ref{lemma:jacobi-functions} is not novel and has already found use in the
literature. In \cite{guo2009} the `generalized Jacobi polynomials/functions' are
denoted $j_n^{(\alpha,\beta)}$, and are defined for all
$\alpha,\beta\in\mathbb{R}$ as 

\begin{align*}
j_n^{(\alpha,\beta)} \propto \left\{
  \begin{array}{lcl}
    P_{n_1}^{(-\alpha,-\beta,-\alpha,-\beta)}, & & \text{$\alpha\leq -1$ and
      $\beta\leq -1$,} \\
    P_{n_1}^{(-\alpha,\beta,-\alpha,0)}, & & \text{$\alpha\leq -1$ and
      $\beta>-1$,} \\
    P_{n_1}^{(\alpha,-\beta,0,-\beta)}, & & \text{$\alpha>-1$ and $\beta\leq
      -1$,} \\
    P_{n_1}^{(\alpha,\beta,0,0)}, & & \text{else,}
  \end{array}\right.
\end{align*}

\noindent where the index $n_1$ is defined as 

\begin{align*}
n_1 = \left\{
  \begin{array}{lcl}
    n - \lfloor -\alpha \rfloor - \lfloor -\beta \rfloor, & & \text{$\alpha\leq -1$ and
      $\beta\leq -1$,} \\
    n - \lfloor -\alpha \rfloor, & & \text{$\alpha\leq -1$ and
      $\beta>-1$,} \\
    n - \lfloor -\beta \rfloor, & & \text{$\alpha>-1$ and $\beta\leq
      -1$,} \\
    n, & & \text{else,}
  \end{array}\right.
\end{align*}

\noindent and the integer floor function is denoted $\lfloor\cdot\rfloor$. These
functions are only defined for certain values of $n$ but \cite{guo2009} presents
significant approximation theory using them. They are advantageous for
solving high-order differential equations with boundary conditions via a global
spectral expansion.

Finally, we present two classical notational conventions that we will use
briefly in the next section. The classical Jacobi polynomials that result from the cases
$\alpha = \beta = - \frac{1}{2}$ and $\alpha = \beta = + \frac{1}{2}$ are the
Chebyshev polynomials of the first and second kinds, respectively. Recalling
the relation $r = \cos \theta$, these polynomials are typically denoted $T_n
(r)$ and $U_n (r)$ and they have a very special and concise representation as
trigonometric polynomials:
\begin{align*}
    \begin{array}{ccccccc}
     \sqrt{\frac{\pi}{2}}  \tilde{P}_n^{(- 1 / 2, - 1 / 2)} (r) &=& T_n (r)
     &=&\cos \left( n \nonesep \theta \right) &=& \cos \left[ n \nonesep
     \arccos (r) \right]\\
     & & &  & & &\\
     \sqrt{\frac{\pi}{2}} \tilde{P}_n^{(1 / 2, 1 / 2)} (r) &=& U_n (r) &=&
     \frac{\sin \left[ (n + 1) \theta \right]}{\sin \theta} &=& \frac{\sin
     \left[ (n + 1) \arccos (r) \right]}{\sin \left[ \arccos \left( r \right)
     \right]} .
     \end{array}
\end{align*}

\subsection{Generalizing the Fourier basis}
\label{sec:derivationFourier}
In this section we will generalize the canonical Fourier basis given by
\[ \Psi_k (\theta) = e^{i \nonesep k \nonesep \theta} . \]
Our methodology is based upon the following dissection of the 
Fourier basis for $k \neq 0$:
\[ \begin{array}{rcccc}
     e^{i \nonesep k \nonesep \theta} & = & \cos \left( k \nonesep \theta)
     \right. & + & i \sin \left( k \nonesep \theta \right)\\
     &  &  &  & \\
     & = & \cos \left( |k| \theta) \right. & + & i \tmop{sgn} (k) \sin \left(
     |k| \theta) \right.\\
     &  &  &  & \\
     & = & T_{|k|} \left( \cos \nonesep \theta \right) & + & i \tmop{sgn} (k)
     \sin (\theta) U_{|k| - 1} \left( \cos \theta \right)\\
     &  &  &  & \\
     & = & \sqrt{\frac{\pi}{2}}  \Big[ \underbrace{\tilde{P}_{|k|}^{(- 1 / 2, - 1 / 2)}
     \left( \cos \theta \right)}_{(a)}  & + & i \tmop{sgn} (k) \underbrace{\sin
     (\theta) \tilde{P}_{|k| - 1}^{\left( 1 / 2, 1 / 2 \right)} \left( \cos
     \theta \right)}_{(b)} \Big] .
   \end{array} \]
We have broken down the Fourier basis into two components: the first component
$(a)$ is even with respect to $\theta$ as it is simply a polynomial in $\cos
\theta$. The second term $(b)$ is odd in $\theta$ as it is a polynomial in
$\cos \theta$ (an even function) multiplied by the odd function $\sin \theta$.
This breakdown suggests that we can construct more general kinds of
Fourier-type functions by augmenting the type of polynomials employed.

However, we cannot switch around polynomials with impunity; we want to
retain orthogonality (at least with respect to some weight function). The
separation into terms $(a)$ and $(b)$ above elucidates the biorthogonal
decomposition of the Fourier basis. The $(a)$ functions are orthogonal with
respect to each other, and with respect to the $(b)$ functions. In this case,
the biorthogonality is manifested as an even-odd separation. Suppose we wish to
generate a basis set orthogonal under the weight $1+ \cos \theta=1+r$. Naturally
we can do this for basis $(a)$ by changing the second Jacobi class parameter
from $\beta = -\frac{1}{2}$ to $\beta=+\frac{1}{2}$. In order to do this for
basis $(b)$, we use Lemma \ref{lemma:jacobi-functions}.

For $\alpha, \beta > - 1$, we have the polynomials $\tilde{P}^{(\alpha,
\beta)}_n$ that are orthogonal in $L^2 \left( [- 1, 1], \mathbbm{R} ;
w_r^{(\alpha, \beta)} \right)$. By setting $a = b = \frac{1}{2}$ in Lemma
\ref{lemma:jacobi-functions}, we also observe that the Jacobi functions
$\tilde{P}^{(\alpha + 1, \beta + 1, 1 / 2, 1 / 2)}_n = (1 - r^2)^{1 / 2} 
\tilde{P}^{(\alpha + 1, \beta + 1)}_n$ are orthogonal under the same
weight. If we set $\alpha = \beta = - \frac{1}{2}$, and add these two functions
together with the appropriate scaling factors, then we exactly recover the
Fourier basis by reversing the dissection steps above (i.e. by creating a
biorthogonal construction). Of course, we are free
to choose any values of $(\alpha, \beta)$ that we desire in order to derive
generalized trigonometric Fourier functions. In fact, this technique has
already been used by Szeg$\ddot{\text{o}}$ {\cite{szego1975}} to determine
orthogonal polynomials on the unit disk. Because the statement in
{\cite{szego1975}} is merely a passing comment and is a markedly different
result than what we desire, we present the following theorem:

\begin{theorem}
  \label{thm:Psifunctions}{\tmem{(cf. Szeg$\ddot{\text{o}}$,
  {\cite{szego1975}})}} For any $\gamma > - \frac{1}{2}$, the functions
  \begin{equation}
    \label{eq:Psi-def} \Psi^{(\gamma)}_k (\theta) = \left\{ \begin{array}{lll}
      \frac{1}{\sqrt{2}} \tilde{P}^{(- 1 / 2, \gamma - 1 / 2)}_0 (\cos
      \theta), &  & k = 0\\
      &  & \\
      \frac{1}{2} \left[ \tilde{P}_{|k|}^{(- 1 / 2, \gamma - 1 / 2)} (\cos
      \theta) + i \tmop{sgn} (k) \sin (\theta) \tilde{P}_{|k| - 1}^{(1 / 2,
      \gamma + 1 / 2)} (\cos \theta) \right], &  & k \neq 0
    \end{array} \right.
  \end{equation}
  are complete and orthonormal in $L^2 \left( [- \pi, \pi], \mathbbm{C} ;
  w_{\theta}^{(\gamma, 0)} \right)$.
\end{theorem}

\begin{proof}
  For orthonormality, it suffices to show
  \begin{enumerate}
    \item $\left\langle \tilde{P}_{|k|}^{(- 1 / 2, \gamma - 1 / 2)} \left(
    \cos \theta \right), \tilde{P}_{|l|}^{(- 1 / 2, \gamma - 1 / 2)} \left(
    \cos \theta \right) \right\rangle_{w_{\theta}^{(\gamma, 0)}} = 2
    \delta_{|k|, |l|}$
    
    \item $\left\langle \sin \theta \tilde{P}_{|k| - 1}^{(1 / 2, \gamma + 1 /
    2)} \left( \cos \theta \right), \sin \theta \tilde{P}_{|l| - 1}^{(1 / 2,
    \gamma + 1 / 2)} \left( \cos \theta \right)
    \right\rangle_{w_{\theta}^{(\gamma, 0)}} = 2 \delta_{|k|, |l|}$, for $k, l
    \neq 0$.
    
    \item $\left\langle \tilde{P}_|k|^{(- 1 / 2, \gamma - 1 / 2)} \left( \cos
    \theta \right), \sin \theta \tilde{P}_{|l| - 1}^{(1 / 2, \gamma + 1 / 2)}
    (\cos \theta) \right\rangle_{w_{\theta}^{(\gamma, 0)}} = 0$, for $l \neq
    0$.
  \end{enumerate}

  The first property is a direct result of orthonormality of the normalized
  Jacobi polynomials $\tilde{P}$ and the observation that on $[0, \pi]$,
  $\left\langle f (\cos \theta), g (\cos \theta)
  \right\rangle_{w_{\theta}^{(\gamma, 0)}} = \left\langle f (r), g (r)
  \right\rangle_{w_r^{(- 1 / 2, \gamma - 1 / 2)}}$. The second property is a
  result of the same observations as the first property along with the result
  of Lemma \ref{lemma:jacobi-functions}. The third property results from the
  fact that an odd function integrated over a symmetric interval is 0. 
  Orthonormality then follows from an explicit calculation of $\left\langle
  \Psi_k^{(\gamma)}, \Psi_l^{(\gamma)} \right\rangle_{w_{\theta}^{(\gamma,
  0)}}$ using the above three properties.

  For completeness we note that any function $f \in L^2$ can be decomposed
  into an even $f_e$ and an odd $f_o$ part. That $f_e$ is representable is
  clear from the fact that $\tilde{P}_n^{(- 1 / 2, \gamma - 1 / 2)} \left(
  \cos \theta \right)$ is complete over $\theta \in [0, \pi]$, which by symmetry
  implies completeness over all $L^2$-even functions $f_e$. Similary, the
  collection of functions $\sin \theta \tilde{P}_n^{(- 1 / 2, \gamma - 1 / 2)}$
  is complete over all $L^2$-odd functions $f_o$ by Lemma
  \ref{lemma:jacobi-functions}. Linearity and orthogonality of the even and odd
  parts yields the result.  \end{proof}

\begin{remark}
  Szeg$\ddot{\text{o}}$ {\cite{szego1975}} gives a more general result that
  involves orthogonality over the weight $w_{\theta}^{(\gamma, \delta)}$ for
  $\delta \neq 0$. We do not require this level of generality; for $\delta
  \neq 0$ the weight function becomes zero at $\theta = 0$, which we will see
  does not help our cause. Indeed, it is possible to generalize
  Szeg$\ddot{\text{o}}$'s result: he derived polynomials on the unit disk
  orthogonal with respect to $w_{\theta}^{(\gamma, \delta)}$. By using Lemma
  \ref{lemma:jacobi-functions} with $a, b$ different from $\frac{1}{2}$, we
  can derive non-polynomial basis sets that are orthogonal under a
  great variety of weights. These functions naturally may not be periodic on
  $\theta \in [- \pi, \pi]$ if the quantity $(1 - r)^a (1 + r)^b$ cannot be
  periodically extended in $\theta$-space to $[- \pi, \pi]$.
\end{remark}

We will refer to the functions (\ref{eq:Psi-def}) as either the generalized
Fourier series, or the Szeg\"{o}-Fourier functions. In the definition of the functions $\Psi_k^{(\gamma)}$ it is desirable to use
the $L^2$-normalized versions of the Jacobi polynomials $\tilde{P}$, rather
than the standard polynomials $P$. If the standard polynomials are used,
then the norm of the Szeg$\ddot{\text{o}}$-Fourier functions
$\Psi_k^{(\gamma)}$ depends on the rather unpleasant-looking sum $h_{|k|}^{(-
1 / 2, \gamma - 1 / 2)} + h_{|k| - 1}^{(1 / 2, \gamma + 1 / 2)}$, and using
this convention implies that $\Psi_k^{(\gamma)}$ is not orthogonal to
$\Psi_{- k}^{(\gamma)}$.

We can also distribute the weight function onto the basis functions, which
gives us orthogonality in the unweighted $L^2$-norm:

\begin{corollary}
  \label{cor:psifunctions}For any $\gamma > - \frac{1}{2}$, the functions
  \[ \psi^{(\gamma)}_k (\theta) = \left\{ \begin{array}{lll}
       \frac{\sqrt[\asterisk]{w_{\theta}^{(\gamma, 0)}}}{\sqrt{2}}
       \tilde{P}^{(- 1 / 2, \gamma - 1 / 2)}_0 (\cos \theta), &  & k = 0\\
       &  & \\
       \frac{\sqrt[\asterisk]{w_{\theta}^{(\gamma, 0)}}}{2} \left[
       \tilde{P}_{|k|}^{(- 1 / 2, \gamma - 1 / 2)} (\cos \theta) + i
       \tmop{sgn} (k) \sin (\theta) \tilde{P}_{|k| - 1}^{(1 / 2, \gamma + 1 /
       2)} (\cos \theta) \right], &  & k \neq 0
     \end{array} \right. \]
  are complete and orthonormal in $L^2 \left( [- \pi, \pi], \mathbbm{C}
  \right)$.
\end{corollary}

\begin{figure}[tbp]
  \begin{center}
  \includegraphics[width=1.0\textwidth]{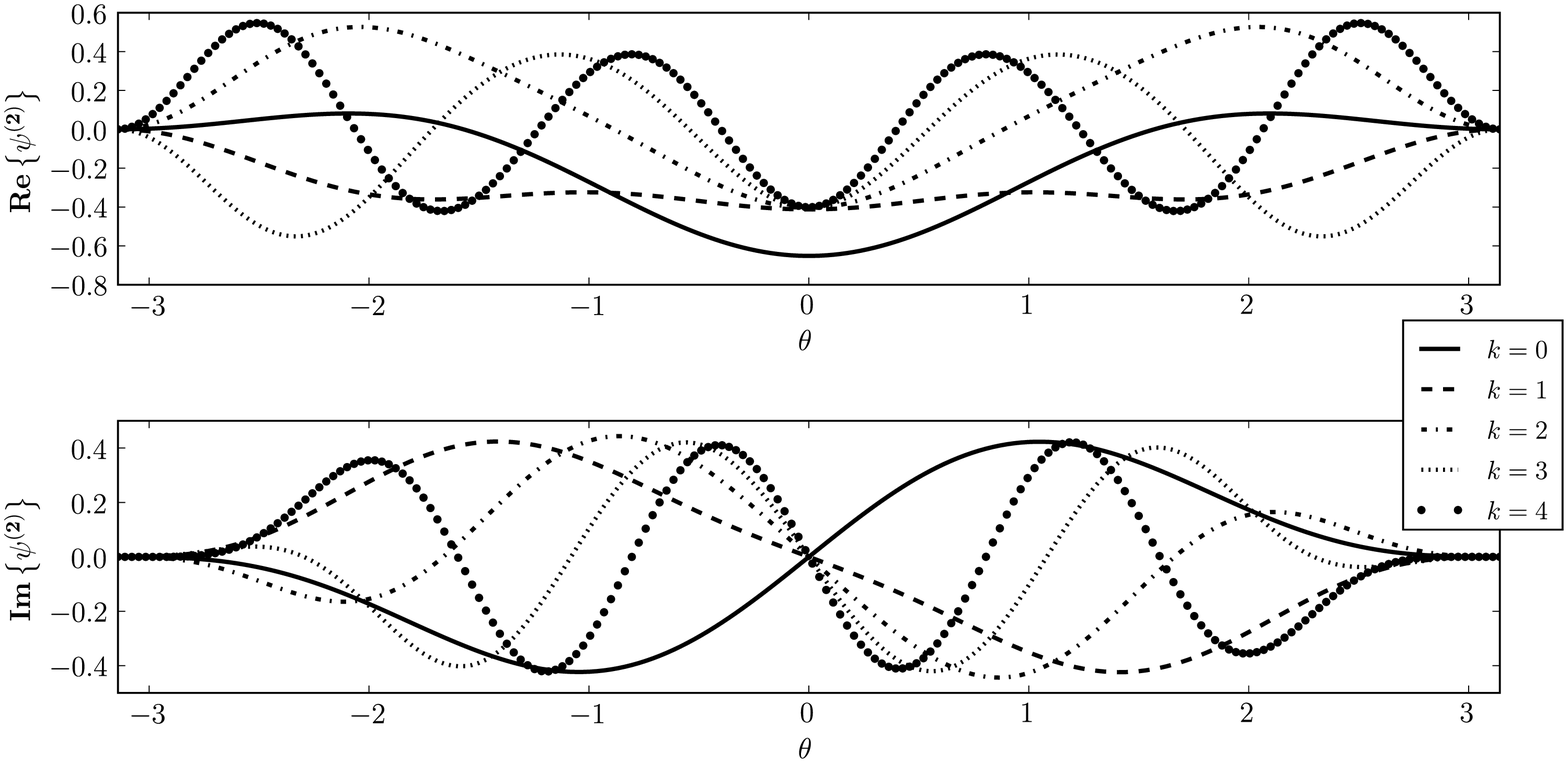}
  \caption{Plots of the weighted Szeg$\ddot{\text{o}}$-Fourier functions
  $\psi^{(2)}_k (\theta)$ for $k = 0, 1, 2, 3,$and $4$. Real part (top) and
  imaginary part (bottom).}
  \label{fig:szegofourier-plots}
  \end{center}
\end{figure}

Due to the properties of $\sqrt[\asterisk]{w_{\theta}^{(\gamma, 0)}}$ given in
(\ref{eq:wt-npbranch}), the functions $\psi^{(\gamma)}_k (\theta)$
decay like $\left( \cos \frac{\theta}{2} \right)^{\gamma}$ at $\theta = \pm
\pi$. This is exemplified in Figure \ref{fig:szegofourier-plots} where we plot
the real and imaginary parts of the functions for $\gamma = 2$. The even/odd
behavior in $\theta$ for real/imaginary components depicted in the figure
depends on the even/odd parity of $\gamma$. (There is no such characterization
possible when $\gamma \nin \mathbbm{N}_0$.) Clearly for $\gamma = 0$ we have
$\Psi_k^{(0)} = \psi_k^{(0)} = \frac{1}{\sqrt{2 \pi}} e^{i \nonesep k \nonesep
\theta}$, the canonical Fourier basis.

\subsection{\label{sec:derivation-map}Mapping to the real line}

Having developed the necessary preliminaries on the finite interval, we now
jump to the infinte line $x \in \mathbbm{R}$ using the mappings introduced in
Table \ref{tab:xforms}. To facilitate the mapping, the following identities
characterizing the mapping between $\theta$-space and $x$-space are useful:
\begin{align*}
     \cos \theta = & \frac{1 - x^2}{1 + x^2}, & (1 - \cos
     \theta) =& \frac{2 x^2}{x^2 + 1},\\
     &  &  & \\
     \sin \theta = & \frac{2 x}{x^2 + 1}, & (1 + \cos \theta) =&
     \frac{2}{x^2 + 1} .
\end{align*}
Using these identities, we rewrite and relabel the functions
$\Psi_k^{(\gamma)} (\theta)$:
\[ \begin{array}{lll}
     \Phi_k^{(s)} (x) & \assign & \Psi_k^{(s - 1)} (\theta (x))\\
     &  & \\
     & = & \left\{ \begin{array}{lll}
       \frac{1}{\sqrt{2}} \tilde{P}^{(- 1 / 2, s - 3 / 2)}_0 \left( \frac{1 -
       x^2}{1 + x^2} \right), &  & k = 0\\
       &  & \\
       \frac{1}{2} \left[ \tilde{P}_{|k|}^{(- 1 / 2, s - 3 / 2)} \left(
       \frac{1 - x^2}{1 + x^2} \right) + \frac{2 i \nonesep x \tmop{sgn}
       (k)}{x^2 + 1} \tilde{P}_{|k| - 1}^{(1 / 2, s - 1 / 2)} \left( \frac{1 -
       x^2}{1 + x^2} \right) \right], &  & k \neq 0
     \end{array} \right.
   \end{array} \]
The above definition is valid for any $s > \frac{1}{2}$. $s = 1$ corresponds
to a mapping of the canonical Fourier basis (i.e., $s\doteq \gamma+1$). These functions are orthogonal
over the weight $w_x^{(s, 0)}$. By following the route from Corollary
\ref{cor:psifunctions} we can distribute the weight over the basis functions,
and in this particular instance we choose the phase-shifted square root given
in (\ref{eq:wx-npbranch}):
\begin{align}
     \phi_k^{(s)} & \assign \sqrt[\asterisk]{w_x^{(s, 0)}} \Phi_k^{(s)}
     (x)\nonumber \\
     & \nonumber \\
     \label{eq:phi-def}
     & = \left\{ \begin{array}{lll}
       \frac{2^{\left( \frac{s - 1}{2} \right)}}{(x - i)^s} \tilde{P}^{(- 1 /
       2, s - 3 / 2)}_0 \left( \frac{1 - x^2}{1 + x^2} \right), &  & k = 0\\
       &  & \\
       \frac{2^{\left( \frac{s}{2} - 1 \right)}}{(x - i)^s} \left[
       \tilde{P}_{|k|}^{(- 1 / 2, s - 3 / 2)} \left( \frac{1 - x^2}{1 + x^2}
       \right) + \frac{2 i \nonesep x \tmop{sgn} (k)}{x^2 + 1} \tilde{P}_{|k|
       - 1}^{(1 / 2, s - 1 / 2)} \left( \frac{1 - x^2}{1 + x^2} \right)
       \right], &  & k \neq 0.
     \end{array} \right.
\end{align}
The functions (\ref{eq:phi-def}) are what we call the generalized Wiener
rational functions. At present there is no clear reason why we have chosen to use
$\sqrt[\asterisk]{w_x^{(s, 0)}}$ instead of the usual square root
$\sqrt{w_x^{(s, 0)}}$ to distribute the weight. However, the corollary
following the coming proposition should provide part of the motivation.

\begin{proposition}
  For any $s > \frac{1}{2}$, the functions $\Phi_k^{(s)} (x)$ are complete and
  orthonormal in $L^2 \left( \mathbbm{R}, \mathbbm{C} ; w_x^{(s, 0)} \right)$.
  The functions $\phi_k^{(s)} (x)$ are complete and orthonormal in $L^2 \left(
  \mathbbm{R}, \mathbbm{C} \right)$. Furthermore, the decay rate of these
  functions can be characterized as
  \[ \lim_{|x| \rightarrow \infty} \left| x^t \phi_k^{(s)} (x) \right| <
     \infty, \hspace{1cm} t \leq s \]
\end{proposition}

\begin{corollary}
  \label{cor:wienerconnection}Recalling the definition of Wiener's original
  basis functions $\phi_n (x)$ in (\ref{eq:wiener-functions}), the
  following relation holds:
  \[ i \sqrt{2} \phi_n^{(1)} (x) \equiv \phi_n (x), \hspace{1cm} n \in
     \mathbbm{N}_0 . \]
\end{corollary}

\begin{figure}[tbp]
  \begin{center}
  \includegraphics[width=1.0\textwidth]{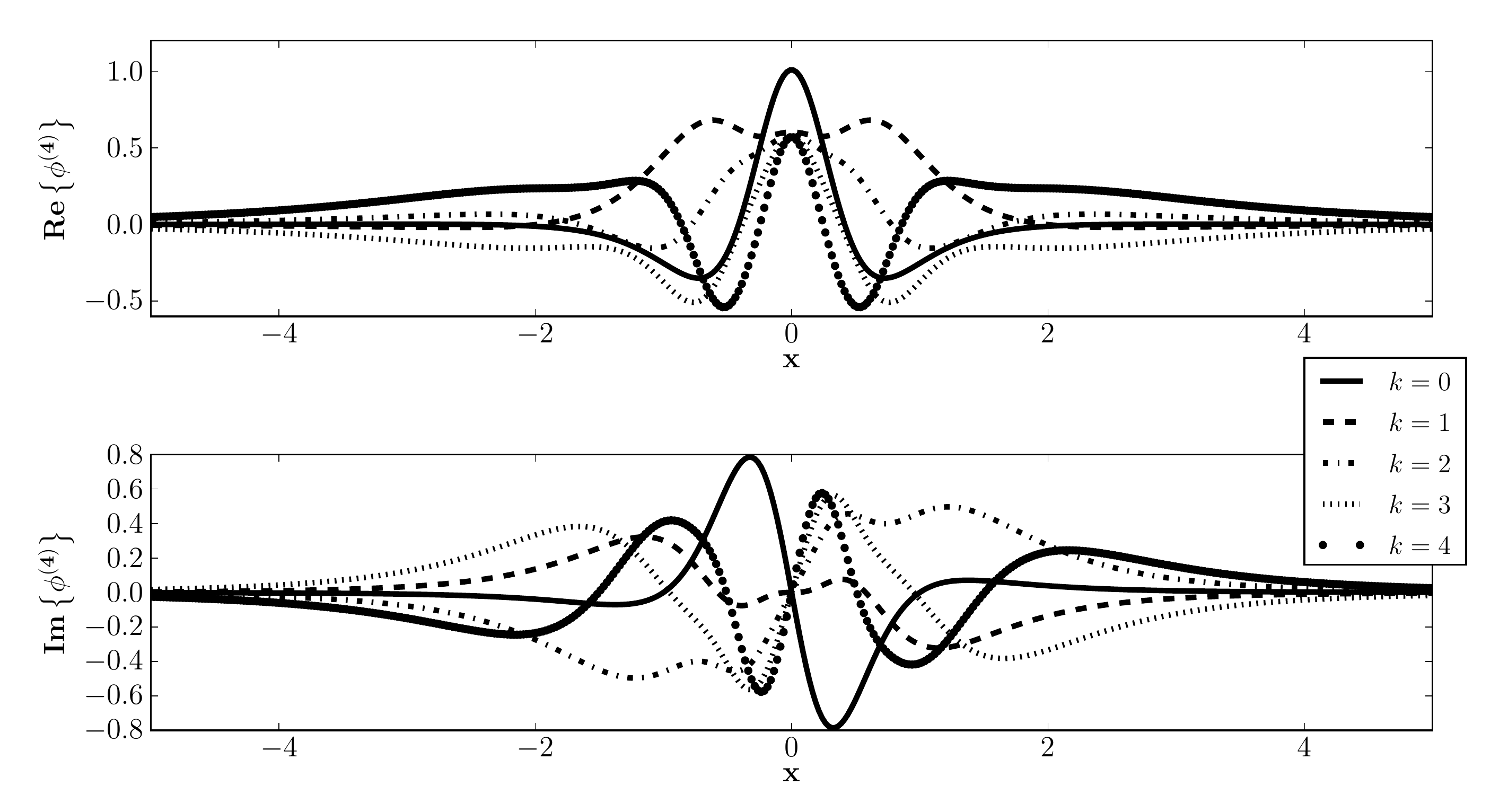}
  \caption{Plots of the functions $\phi_k^{(4)} (x)$ for $k = 0, 1, 2, 3, 4$.}
  \label{fig:wienerplots}
  \end{center}
\end{figure}

\noindent We show plots of the functions $\phi_k^{(4)}$ in Figure \ref{fig:wienerplots}.
The conclusion of the corollary is easily seen if one makes the connection
\[ e^{i \theta} = \frac{i - x}{i + x}, \]
along with knowledge of the fact that $\Phi_k^{(1)} (x) = \psi_k^{(0)}
(\theta) = \frac{1}{\sqrt{2 \pi}} e^{i \nonesep k \nonesep \theta}$. We have
thus shown that the orthogonal functions $\phi_k^{(s)}$ over the real line are
a generalization of Wiener's original basis set. Furthermore, $\phi_k^{(s)}$
decays like $x^{- s}$ while retaining orthogonality under the same unit weight
measure. When $s$ is an integer, the functions are also purely rational: they
are the division of one complex-valued polynomial in $x$ by another. This
connection was rather helpful in the nascent stages of the computing when the
calculation of a non-polynomial function required significantly more
computational investment, but now this property is probably more aesthetic
than functional. As a result, our use of the quantity
$\sqrt[\asterisk]{w_x^{(s, 0)}}$ is not entirely necessary for purposes of
evaluating the functions; it is equally valid to use the traditional squre
root $\sqrt{w_x^{(s, 0)}}$.

By using the traditional square root, one sacrifice made is that the
analogous written form of Corollary \ref{cor:wienerconnection} becomes less
fortuitous and is complicated by $x$-dependent phase-shift factors. The same
observation is true of the weight $\sqrt[\asterisk]{w_{\theta}^{(\gamma, 0)}}$
used in the definition of the Szeg$\ddot{\text{o}}$-Fourier functions
$\psi_k^{(\gamma)} (\theta)$ in Corollary \ref{cor:psifunctions}. A second
reason to use the phase-shifted square root is that it can be written in the
following convenient form:
\begin{equation}
  \label{eq:phase-shift-rewrite} \sqrt[\ast]{w_x^{(s, 0)}} = \left[
  \frac{i}{\sqrt{2}} \left( 1 + e^{- i \theta} \right) \right]^s .
\end{equation}
The utility of this expression will become clear when we consider the
connection problems in Section \ref{sec:propjacobi}.

We have accomplished our goal of deriving basis functions satisfying tunable
decay rate while maintaining $L^2$-orthogonality. However, it is not clear
that these are superior or useful functions.  We will now present some
properties of the basis and make the argument that these basis functions indeed
are very useful for solving problems in scientific computing.

%% file: part1Content/fourierprop.tex
\label{sec:prop}In this section we explore some of the desirable properties of
the generalized Wiener basis set $\left\{ \phi_k^{(s)} \right\}_{k \in
\mathbbm{Z}}$, $s > \frac{1}{2}$ based on their close relation to Fourier
Series. The argument we make is that these functions inherit all the useful
properties of the Fourier basis with the additional property that the decay
rate $s$ at $|x| = \infty$ may be chosen. Many of these
properties (e.g. the sparse modal differentiation matrix) rely on Jacobi
polynomial properties covered in the next section. 
In particular, although application of the FFT is indeed a virtue of this basis, we will
discuss it only in Part \Rmnum{2}, which focuses with computational issues.

\subsection{Symmetry}

The derivation of the basis functions automatically yields various simple
properties. Note that due to the mapping, any property of the basis on the
real line $x \in \mathbbm{R}$ also applies to the respective trigonometric
interval $\theta \in [- \pi, \pi]$. We omit the proof of these properties as
they are elementary:
\begin{enumerate}
  \item Index symmetry
  \begin{eqnarray}
    \Phi_k^{(s)} (x) & = & \overline{\Phi_{- k}^{(s)} (x)} 
    \label{eq:indsym-1}\\ 
    &  &  \nonumber\\
    | \Phi_k^{(s)} (x) | & = & | \Phi_{- k}^{(s)} (x) | \nonumber\\
    &  &  \nonumber\\
    | \phi_k^{(s)} (x) | & = & | \phi_{- k}^{(s)} (x) | \nonumber\\
    &  &  \nonumber\\
    \phi_k^{(1)} (x) & = & \overline{\phi_{- k - 1}^{(1)} (x)} \nonumber
  \end{eqnarray}
  \item Function symmetry
  \[ \begin{array}{ccc}
       \tmop{Re} \left\{ \Phi_k^{(s)} (x) \right\} & = & \tmop{Re} \left\{
       \Phi_k^{(s)} (- x) \right\}\\
       &  & \\
       \tmop{Im} \left\{ \Phi_k^{(s)} (x) \right\} & = & - \tmop{Im} \left\{
       \Phi_k^{(s)} (- x) \right\}\\
       &  & \\
       | \Phi_k^{(s)} (x) | & = & | \Phi_k^{(s)} (- x) |\\
       &  & \\
       | \phi_k^{(s)} (x) | & = & | \phi_k^{(s)} (- x) |
     \end{array} \]
\end{enumerate}

\subsection{Periodicity}

Because trigonometric polynomials are periodic over $\theta \in [- \pi, \pi]$,
we cannot expect this condition to be violated on the infinite interval $x \in
\mathbbm{R}$. From the viewpoint of expanding functions over the infinite
interval $\mathbbm{R}$, the points $x = \pm \infty$ are both unique points.
However, because of the mapping, the basis functions view the points $x = \pm
\infty$ the same as they view the points $\theta = \pm \pi$: i.e. they are the
same point. This serves as a disadvantage if we wish to e.g. expand functions
with different decay rates at $\pm \infty$ because this is in effect
non-smooth behavior of the function at a single point, which degrades the
convergence rate of the approximation.

In particular it is known that although a Fourier series approximation will
converge in the $L^2$ sense for an $L^2$ function, the rate of convergence is
only algebraic if the function is non-periodic. Naturally, this deficiency will follow us to the infinite
interval. Indeed, such observations have already been made
{\cite{canuto2006}}. Empirical studies we have carried out show that the
concern of periodicity is not paramount and frequently one can overlook it when
comparing results to other basis expansions. Nonperiodic behavior is often
manifested as algebraic decay at $x = \pm\infty$, where existing basis sets
already have problems in approximation. In Part \Rmnum{2} we will present
examples that explicitly illustrate this lack of fast convergence rate when the
function to be expanded is not `periodic' at $x = \pm \infty$.

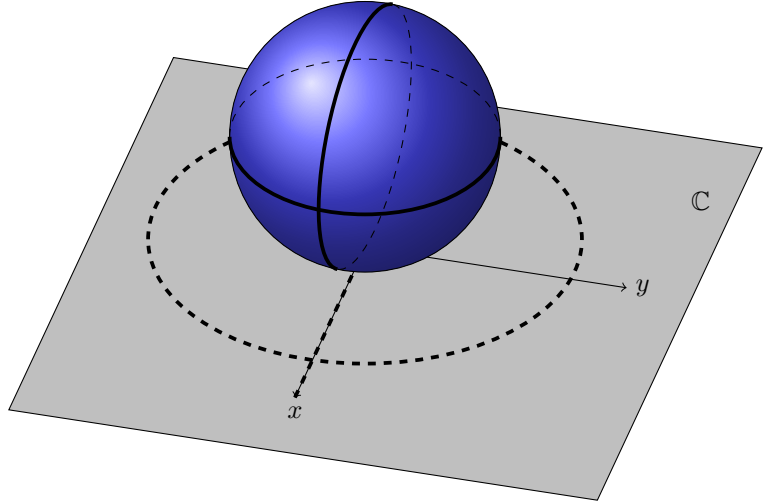
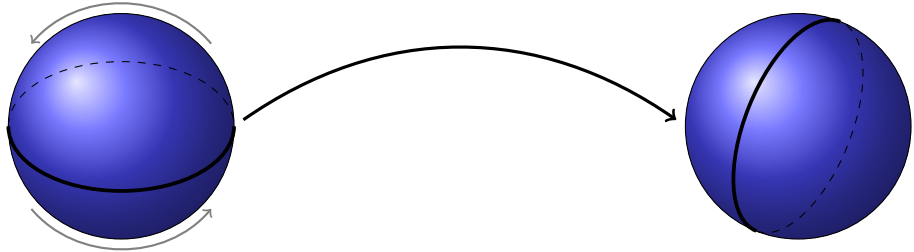
\begin{figure}[htb]
 \begin{center}
  \subfigure[Illustration of the stereographic connection between the Riemann
  Sphere and the complex plane. The equator corresponds to the unit circle, and
  the meridian can be identified with the real ($x$) axis.]{
  \input{tikz/part1/rsphere}\label{fig:rsphere}}

  \subfigure[The effect of the linear fraction map we've chosen to take $\theta =
    \arg z_1$
  to $x = {\rm Re} \{z_2\}$: a rotation of the Riemann
  Sphere.]{\input{tikz/part1/riemann_rotate}\label{fig:riemannrotate}}
  \caption{The linear fractional mapping that relates $x$ to $\theta$ has an
  illuminating representation when viewed as a transformation of the complex
  plane.}
  \end{center}
\end{figure}

Note that although it may seem a bit unnatural that periodicity is a
condition at $x = \pm \infty$, in fact it is not surprising at all. One
may consider our mapping as a rather unremarkable tangent map from
$\theta$-space to $x$-space as written in Table \ref{tab:xforms}. However, it is
more deep than that: the functions $\Psi_k^{(\gamma)} (\theta)$ and
$\psi_k^{(\gamma)} (\theta)$ are periodic basis sets for complex-valued
functions on the unit circle $\mathbbm{T}$. In other words, we may actually view
these basis sets as functions of $z \in \mathbbm{C}$. What looks like a tangent
mapping from $\theta$-space to $x$-space is actually a linear fractional map (a
M\"{o}bius transformation) from the unit circle (the complex plane) to
complexified $x$-space (the complex plane).

Linear fractional maps are structure-preserving maps of the complex plane;
an illuminating way to consider this is by identifying the complex plane
$\mathbbm{C}$ with the Riemann Sphere (see Figure \ref{fig:rsphere}). Then
the linear fractional map we've chosen merely takes one great circle
(the unit circle in $z$-space) to another great circle (the real line
in complexified $x$-space). Therefore, our approximation is nothing
more than a rotation of functions on the Riemann Sphere (see Figure
\ref{fig:riemannrotate}); the target space simply happens to correspond
to the real line. From this point of view, periodicity at $|x| = \infty$
(analyticity at complexified $x = \infty$) is natural.

Nevertheless, this `natural' periodicity can be problematic if we attempt to
approximate a function that is not complex-analytic at $x = \infty$. In Part
\Rmnum{2} we present examples of functions that are not analytic at $x = \infty$ and
we will empirically analyze the impact of violating the assumption of
periodicity.


%% file: tikz/part1/rsphere.tex
\newcommand\pgfmathsinandcos[3]{%
  \pgfmathsetmacro#1{sin(#3)}%
  \pgfmathsetmacro#2{cos(#3)}%
}
\newcommand\LongitudePlane[3][current plane]{%
  \pgfmathsinandcos\sinEl\cosEl{#2} 
  \pgfmathsinandcos\sint\cost{#3} 
  \tikzset{#1/.estyle={cm={\cost,\sint*\sinEl,0,\cosEl,(0,0)}}}
}
\newcommand\LatitudePlane[3][current plane]{%
  \pgfmathsinandcos\sinEl\cosEl{#2} 
  \pgfmathsinandcos\sint\cost{#3} 
  \pgfmathsetmacro\yshift{\cosEl*\sint}
  \tikzset{#1/.estyle={cm={\cost,0,0,\cost*\sinEl,(0,\yshift)}}} %
}
\newcommand\DrawLongitudeCircle[2][1]{
  \LongitudePlane{\angEl}{#2}
  \tikzset{current plane/.prefix style={scale=#1}}
  \pgfmathsetmacro\angVis{atan(sin(#2)*cos(\angEl)/sin(\angEl))} %
  \draw[current plane,line width=0.5mm] (\angVis:1) arc (\angVis:\angVis+180:1);
  \draw[current plane,dashed] (\angVis-180:1) arc (\angVis-180:\angVis:1);
}
\newcommand\DrawLatitudeCircle[2][1]{
  \LatitudePlane{\angEl}{#2}
  \tikzset{current plane/.prefix style={scale=#1}}
  \pgfmathsetmacro\sinVis{sin(#2)/cos(#2)*sin(\angEl)/cos(\angEl)}
  \pgfmathsetmacro\angVis{asin(min(1,max(\sinVis,-1)))}
  \draw[current plane,line width=0.5mm] (\angVis:1) arc (\angVis:-\angVis-180:1);
  \draw[current plane,dashed] (180-\angVis:1) arc (180-\angVis:\angVis:1);
}

\newcommand\DrawPlanarCircle[2][1]{
  \LatitudePlane{\angEl}{#2}
  \tikzset{current plane/.prefix style={scale=#1}}
  \pgfmathsetmacro\sinVis{sin(#2)/cos(#2)*sin(\angEl)/cos(\angEl)}
  \pgfmathsetmacro\angVis{asin(min(1,max(\sinVis,-1)))}
  \draw[current plane,dashed,line width=0.5mm] (180:1) arc (180:-180:1);
}

\begin{tikzpicture}
\begin{scope}[scale=.90,xshift=1cm,yshift=1cm]

\def\R{2.0} 
\def\angEl{35} 
\def\angAz{-105} 
\def\angPhi{-40} 
\def\angBeta{19} 


\pgfmathsetmacro\H{\R*cos(\angEl)} 
\tikzset{xyplane/.estyle={cm={cos(\angAz),sin(\angAz)*sin(\angEl),-sin(\angAz),
                              cos(\angAz)*sin(\angEl),(0,-\H)}}}

\draw[xyplane,fill=gray!50!white] (-2*\R,-2*\R) rectangle (2.7*\R,2.5*\R);
\draw[xyplane] (-1.0*\R,2.3*\R) node[above] {$\mathbb{C}$};
\draw[xyplane,<->] (2.0*\R,0) node[below] {$x$} -- (0,0) -- (0,2.0*\R)
    node[right] {$y$};
\DrawPlanarCircle[\R+1.7]{-30}
\draw[xyplane,line width=0.5mm, dashed] (2.0*\R,0) -- (0,0);

\filldraw[ball color=blue!70!white] (0,0) circle (\R);

\foreach \t in {0} { \DrawLatitudeCircle[\R]{\t} }
\foreach \t in {-110} { \DrawLongitudeCircle[\R]{\t} }



\end{scope}
\end{tikzpicture}

%% file: tikz/part1/riemann_rotate.tex
\newcommand\pgfmathsinandcos[3]{%
  \pgfmathsetmacro#1{sin(#3)}%
  \pgfmathsetmacro#2{cos(#3)}%
}
\newcommand\LongitudePlane[3][current plane]{%
  \pgfmathsinandcos\sinEl\cosEl{#2} 
  \pgfmathsinandcos\sint\cost{#3} 
  \tikzset{#1/.estyle={cm={\cost,\sint*\sinEl,0,\cosEl,(0,0)}}}
}
\newcommand\LatitudePlane[3][current plane]{%
  \pgfmathsinandcos\sinEl\cosEl{#2} 
  \pgfmathsinandcos\sint\cost{#3} 
  \pgfmathsetmacro\yshift{\cosEl*\sint}
  \tikzset{#1/.estyle={cm={\cost,0,0,\cost*\sinEl,(0,\yshift)}}} %
}
\newcommand\DrawLongitudeCircle[2][1]{
  \LongitudePlane{\angEl}{#2}
  \tikzset{current plane/.prefix style={scale=#1}}
  \pgfmathsetmacro\angVis{atan(sin(#2)*cos(\angEl)/sin(\angEl))} %
  \draw[current plane,line width=0.5mm] (\angVis:1) arc (\angVis:\angVis+180:1);
  \draw[current plane,dashed] (\angVis-180:1) arc (\angVis-180:\angVis:1);
}
\newcommand\DrawLatitudeCircle[2][1]{
  \LatitudePlane{\angEl}{#2}
  \tikzset{current plane/.prefix style={scale=#1}}
  \pgfmathsetmacro\sinVis{sin(#2)/cos(#2)*sin(\angEl)/cos(\angEl)}
  \pgfmathsetmacro\angVis{asin(min(1,max(\sinVis,-1)))}
  \draw[current plane,line width=0.5mm] (\angVis:1) arc (\angVis:-\angVis-180:1);
  \draw[current plane,dashed] (180-\angVis:1) arc (180-\angVis:\angVis:1);
}

\newcommand\DrawPlanarCircle[2][1]{
  \LatitudePlane{\angEl}{#2}
  \tikzset{current plane/.prefix style={scale=#1}}
  \pgfmathsetmacro\sinVis{sin(#2)/cos(#2)*sin(\angEl)/cos(\angEl)}
  \pgfmathsetmacro\angVis{asin(min(1,max(\sinVis,-1)))}
  \draw[current plane,dashed,line width=0.5mm] (180:1) arc (180:-180:1);
}

\begin{tikzpicture}[scale=1.5]
\def\angEl{35} 
\def\angAz{-105} 
\def\angPhi{-40} 
\def\angBeta{19} 
\def\R{1}

\pgfmathsetmacro\H{\R*cos(\angEl)} 
\tikzset{xyplane/.estyle={cm={cos(\angAz),sin(\angAz)*sin(\angEl),-sin(\angAz),
                              cos(\angAz)*sin(\angEl),(0,-\H)}}}
\filldraw[ball color=blue!70!white] (2,-0.5) circle (\R);
\node (ballone) at (2+\R,-0.5) {};
\begin{scope}[xshift=2cm,yshift=-0.5cm]
  \DrawLatitudeCircle[\R]{0}
  \begin{scope}[yshift=0.15cm]
    \node (ur) at (30:\R) {};
    \node (ul) at (150:\R) {};
  \end{scope}
  \begin{scope}[yshift=-0.15cm]
    \node (ll) at (210:\R) {};
    \node (lr) at (330:\R) {};
  \end{scope}
\end{scope}

\draw [->,color=gray,thick] (ur) to [out=130,in=50] (ul);
\draw [->,color=gray,thick] (ll) to [out=310,in=230] (lr);

\filldraw[ball color=blue!70!white] (8,-0.5) circle (\R);
\node (balltwo) at (8-\R,-0.5) {};
\begin{scope}[xshift=8cm,yshift=-0.5cm]
  \DrawLongitudeCircle[\R]{55}
\end{scope}

\draw [->,very thick] (ballone) to [out=35,in=145] (balltwo);

\end{tikzpicture}

%% file: part1Content/jacobiprop.tex
\label{sec:propjacobi}The generalized Wiener functions are composed of Jacobi
polynomials, and so it is reasonable to expect that we can use the properties
of the Jacobi polynomials to perform certain tasks using the Wiener basis.
Indeed, we can form recurrence relations, connection coefficients, a Gauss-like
quadrature, and obtain an extremely useful sparsity result for the Galerkin
stiffness matrix.

\subsection{Recurrence Relations}

Due to the strong dependence of the Szeg$\ddot{\text{o}}$-Fourier functions on
the Jacobi polynomials, they inherit six-term recurrence relations from the
three-term recurrences for orthogonal polynomials.
\[ \begin{array}{lll}
     D_n^{(\gamma)} \Psi_{n + 1}^{(\gamma)} & = & \left[ A_n^{(\gamma)} e^{i
     \theta} - B_n^{(\gamma)} \right] \Psi_n^{(\gamma)} + \left[ A_{-
     n}^{(\gamma)} e^{- i \theta} - B_{- n}^{(\gamma)} \right] \Psi_{-
     n}^{(\gamma)} + C_n^{(\gamma)} \Psi_{n - 1}^{(\gamma)} + C_{-
     n}^{(\gamma)} \Psi_{- (n - 1)}^{(\gamma)},\\
     &  & \\
     \Psi_{n + 1}^{(\gamma)} & = & \left[ U_n^{(\gamma)} \cos \theta -
     V_n^{(\gamma)} \right] \Psi_n^{(\gamma)} + \left[ U_{- n}^{(\gamma)} \cos
     \theta - V_{- n}^{(\gamma)} \right] \Psi_{- n}^{(\gamma)} +
     W_n^{(\gamma)} \Psi_{n - 1}^{(\gamma)} + W_{- n}^{(\gamma)} \Psi_{- (n -
     1)}^{(\gamma)},\\
     &  & \\
     \Psi_{n + 1}^{(\gamma)} & = & \left[ \tilde{U}_n^{(\gamma)} i \sin \theta
     - \tilde{V}_n^{(\gamma)} \right] \Psi_n^{(\gamma)} + \left[ \tilde{U}_{-
     n}^{(\gamma)} i \sin \theta - \tilde{V}_{- n}^{(\gamma)} \right] \Psi_{-
     n}^{(\gamma)} + \tilde{W}_n^{(\gamma)} \Psi_{n - 1}^{(\gamma)} +
     \tilde{W}_{- n}^{(\gamma)} \Psi_{- (n - 1)}^{(\gamma)} .
   \end{array} \]
We give formulae for all the real-valued constants $A, B, C, D, U, V, W,
\tilde{U}, \tilde{V}, \tilde{W}$ in Appendix \ref{app:recurrence}. Note that
since the $\Psi_k^{(\gamma)}$ are \tmtextit{not} polynomials in $z = e^{i
\theta}$, there is not a three-term recurrence as there would normally be for
orthogonal polynomials on the unit disk (unless of course $\gamma = 0$).
Although the above formulae are complex-valued six-term recurrence relations,
they are no more difficult computationally than the pair of three-term
recurrences necessary to generate $\tilde{P}_n^{(\alpha, \beta)}$ and
$\tilde{P}_n^{(\alpha + 1, \beta + 1)}$ because $\Psi_n^{(\gamma)}$ is the
complex conjugate of $\Psi_{- n}^{(\gamma)}$ and therefore does not need to be
generated independently. Direct use of any of the above
six-term recurrences for generating the $\Psi_k^{(\gamma)}$ is just as
expensive as forming $\Psi_k^{(\gamma)}$ by the even/odd synthesis of
$\tilde{P}_n^{(\alpha, \beta)}$ and $\tilde{P}_n^{(\alpha + 1, \beta + 1)}$ in
Theorem \ref{thm:Psifunctions}. However, using presumably existing routines
for evaluating Jacobi polynomials and then synthesizing them is likely easier
from an implementation view.

It is reassuring to note that simplifying the recurrence constants in the case
$\gamma = 0$ yields, up to normalization, the trivial recurrence relations for
the monomials on the unit disk $\Psi_k^{(0)} (\arg z) = \frac{z^k}{\sqrt{2 \pi}}$:
\[ \begin{array}{lll}
     \Psi_{n + 1}^{(0)} & = & e^{i \theta} \Psi_n^{(0)},\\
     &  & \\
     \Psi_{n + 1}^{(0)} & = & 2 \cos \theta \Psi_n^{(0)} - \Psi_{n -
     1}^{(0)},\\
     &  & \\
     \Psi_{n + 1}^{(0)} & = & 2 i \sin \theta \Psi_n^{(0)} + \Psi_{n -
     1}^{(0)} .
   \end{array} \]

\noindent Naturally, a recurrence relation for the unweighted $\Psi_k^{(\gamma)}
(\theta)$ translates directly into one for the unweighted Wiener rational
functions $\Phi_k^{(s)} (x)$. The weighted functions $\psi_k^{(\gamma)}
(\theta)$ and $\phi_k^{(s)} (x)$ can be generated by first generating the
unweighted functions and then multiplying by the phase-shifted square root
$\sqrt[\asterisk]{w}$.

\subsection{Connection Problems}

\label{sec:propjacobi-connections}One advantage in using the generalized
Wiener rational function basis is the ability to choose the parameter $s$,
which indicates the rate of decay. In many applications, it may be useful to
augment the basis functions mid-computation to suit the dynamics occuring at a
particular time. In this case, one would like to be able to transfer from one
basis to another while keeping the (finite-term) function expansion identical.
We will also see in Part \Rmnum{2} that this problem also appears in an algorithm
utilizing the FFT. In classical orthogonal polynomial theory, the problem of
equating one expansion to another boils down to determining the connection
coefficients. Before undertaking this task, we first outline the major
tasks we wish to perform.

There are two main tasks on the infinite interval that require connections of
some form:
\begin{enumerate}
  \item Usage of the fast Fourier transform -- transforming $N$ nodal
  evaluations into $N$ modal coefficients (or vice-versa) for an expansion in $\phi^{(s)}$.
  
  \item For a given expansion in $\phi^{(s)}$ (i.e. a set of modal
  coefficients), translating this into a modal coefficient expansion in
  $\phi^{(S)}$ for some $s \neq S$.
\end{enumerate}
In Part \Rmnum{2} where we outline computational considerations, we will address the
above tasks. However, for now it suffices to note that these two tasks can be
reduced to the following three connection problems in $\theta$-space:
\begin{enumerate}
  \item The $\Psi^{(\gamma)}$-$\Psi^{(\Gamma)}$ connection (a necessary
  ingredient for all connection-like tasks)
  
  \item The $\Psi^{(\gamma)}$-$\psi^{(\gamma)}$ connection (a generalization
  of the FFT task)
  
  \item The $\psi^{(\gamma)}$-$\psi^{(\Gamma)}$ connection (identical to
  modification of $s$)
\end{enumerate}
In Sections
\ref{sec:propjacobi-connections-fourier}-\ref{sec:propjacobi-connections-smod},
we will tackle each of these problems. Note that modification of any of the
following finite-interval algorithms for the infinite interval is trivial: the
relations $\Psi_k^{(\gamma)} (\theta) \equiv \Phi_k^{(s - 1)} (x)$,
$\psi_k^{(\gamma)} (\theta) \equiv \phi_k^{(s - 1)} (x)$, and $\gamma \assign
s - 1$ allows for us to easily employ the same operations, whether we want to
do it in $\theta$-space or $x$-space.

\subsubsection{The $\Psi$-$\Psi$ Connection Problem}

\label{sec:propjacobi-connections-fourier}Suppose we have a function $f \in
L^2 \left( [- \pi, \pi], \mathbbm{C} ; w_{\theta}^{(\gamma)} \right) \bigcap
L^2 \left( [- \pi, \pi], \mathbbm{C} ; w_{\theta}^{(\Gamma)} \right)$ with a
Fourier expansion for some $\gamma > - \frac{1}{2}$:
\[ f (x) = \sum_{k \in \mathbbm{Z}} \hat{f}^{(\gamma)}_k \Psi_k^{(\gamma)} .
\]
The goal is determine a way to re-expand $f$ in a Fourier expansion for a
different decay parameter $\Gamma$:
\[ f (x) = \sum_{k \in \mathbbm{Z}} \hat{f}^{(\Gamma)}_k \Psi_k^{(\Gamma)} .
\]
The shift $\Gamma - \gamma$ can take values in the interval $\left( -
\frac{1}{2} - \gamma, \infty \right)$. Naturally one may equate the two
expansions and use orthogonality to relate one set of expansion coefficients
to the other:
\[ \hat{f}^{(\Gamma)}_k = \sum_{l \in \mathbbm{Z}} \hat{f}^{(\gamma)}_l
   \left\langle \Psi_l^{(\gamma)}, \Psi_k^{(\Gamma)}
   \right\rangle_{w_{\theta}^{(\Gamma)}} . \]
We can then define the connection coefficients
\[ \lambda^{\Psi}_{k, l} = \left\langle \Psi_l^{(\gamma)}, \Psi_k^{(\Gamma)}
   \right\rangle_{w_{\theta}^{(\Gamma)}}, \]
where we have suppressed the dependence of $\lambda$ on $\gamma$ and $\Gamma$.
Our task is to determine how to calculate these connection coefficients. Due
to orthogonality, it is clear that
\begin{equation}
  \label{eq:Psi-connection-orthogonality} \lambda_{k, l}^{\Psi} \equiv 0,
  \hspace{1cm} |l| < |k|.
\end{equation}
This implies that the connection problem is solved via the relation
\begin{equation}
  \label{eq:Psi-connection-infinite} \hat{f}^{(\Gamma)}_k =
  \sum_{\tmscript{\begin{array}{c}
    l \in \mathbbm{Z},\\
    |l| \geq |k|
  \end{array}}} \hat{f}^{(\gamma)}_l \lambda_{k, l}^{\Psi} .
\end{equation}
Relation (\ref{eq:Psi-connection-infinite}) is still not attractive: we must
perform an infinite number of operations for an exact connection. If we only
have a finite expansion (say a total of $N$ modal coefficients), we must still
perform $\mathcal{O} (N^2)$ operations to capture all the information at our disposal.
However we will show that, for integer values of the shift $\Gamma - \gamma$,
the connection problem can be solved inexpensively. To be precise, we will
show that for $G \in \mathbbm{N}$, (\ref{eq:Psi-connection-infinite}) reduces
to
\begin{equation}
  \label{eq:Psi-connection-sparse} \hat{f}_k^{(\gamma + G)} = \sum_{k + G \geq
  |l| \geq |k|} \hat{f}^{(\gamma)}_l \lambda_{k, l}^{\Psi} .
\end{equation}
That is, only $2 (G + 1)$ operations per coefficient are necessary to solve
the connection problem (independent of $k$, and of any truncation size $N$).
We refer to the above collapse of the infinite connection problem
(\ref{eq:Psi-connection-infinite}) into the finite $N$-indepedent problem
(\ref{eq:Psi-connection-sparse}) as a sparse connection.

In order to relate one Fourier function to another, we first recall a result
from {\cite{narayan2009}} using (\ref{eq:jdemotiona}) -- (\ref{eq:jpromotionb})
that states that the connection coefficients binding one Jacobi polynomial class
to another are sparse in certain special circumstances.

\begin{lemma}
  \label{lemma:jacobi-connection}For any $\alpha, \beta > - 1$ and any $A, B,
  \in \mathbbm{N_0}_0$, the connection problem
  \[ f (r) = \sum_{n = 0}^{\infty} \hat{f}^{(\alpha, \beta)}_n 
     \tilde{P}^{(\alpha, \beta)}_n (r) \hspace{0.3cm} \longrightarrow
     \hspace{0.3cm} f (r) = \sum_{n = 0}^{\infty} \hat{f}^{(\alpha + A, \beta
     + B)}_n \tilde{P}^{(\alpha + A, \beta + B)}_n (r), \]
  can be solved exactly via the relation
  \begin{equation}
    \label{eq:fhat-connection-sparse} \hat{f}^{(\alpha + A, \beta + B)}_n =
    \sum_{m = 0}^{A + B} \lambda^P_{n, n + m}  \hat{f}_{n + m}^{(\alpha,
    \beta)} .
  \end{equation}
\end{lemma}
In the above we have suppressed the dependence of $\lambda^P$ on $\alpha,
\beta, A,$and $B$, but in the sequel we shall occasionally refer to the above
coefficients as $\lambda^P_{n, m} (\alpha, \beta, A, B)$. The result
(\ref{eq:fhat-connection-sparse}) is not a trivial one; the upper limit for
the sum on the right-hand side is $\infty$ for a general connection problem.
For the very special cases satisfying the lemma, the exact connection becomes
finite. We have not shown how to obtain the Jacobi-Jacobi connection
coefficients $\lambda^P$. For this, one may use explicit formulae given in
{\cite{maroni2008}} or {\cite{askey1975}}, or one may utilize the algorithm
given in {\cite{narayan2009}}.

The above result can be expanded to apply to the
Szeg$\ddot{\text{o}}$-Fourier functions $\Psi_k^{(\gamma)} (\theta)$ and the
corresponding mapped functions $\Phi_k^{(s)} (x)$.

\begin{proposition}
  \label{prop:Psi-connection}For any $\gamma > - \frac{1}{2}$ and any $G \in
  \mathbbm{N}$, the connection problem
  \[ f (\theta) = \sum_{k = - \infty}^{\infty} \hat{f}_k^{(\gamma)}
     \Psi_k^{(\gamma)} (\theta) \hspace{0.3cm} \longrightarrow \hspace{0.3cm}
     f (\theta) = \sum_{k = - \infty}^{\infty} \hat{f}_k^{(\gamma + G)}
     \Psi_k^{(\gamma + G)} (\theta), \]
  can be solved exactly via the relation
  \begin{equation}
    \label{eq:fourier-mode-connection} \hat{f}_k^{(\gamma + G)} = \sum_{l =
    |k|}^{|k| + G} \lambda^{\Psi}_{k, l}  \hat{f}_l^{(\gamma)} + \sum_{l = -
    |k| - G}^{- |k|} \lambda^{\Psi}_{k, l}  \hat{f}_l^{(\gamma)} .
  \end{equation}
\end{proposition}

\noindent Note that (\ref{eq:fourier-mode-connection}) is exactly
(\ref{eq:Psi-connection-sparse}). By making the connection $s - 1
\longleftrightarrow \gamma$, we recover $\lambda^{\Phi}_{k, l} \equiv
\lambda_{k.l}^{\Psi}$, where $\Phi_k^{(s)} (x)$ are the maps of the
Szeg$\ddot{\text{o}}$-Fourier functions $\Psi_k^{(\gamma)}$. We stress again
that this result is nontrivial. This also yields the functional connection
\begin{equation}
  \label{eq:fourier-fun-connection} \begin{array}{lll}
    \Psi_m^{(\gamma)} (\theta) & = & \left\{ \begin{array}{lll}
      \sum_{|k| \leq m} \lambda_{k, m}^{\Psi} \Psi_k^{(\gamma + G)} (\theta),
      &  & |m| \leq G\\
      &  & \\
      \sum_{m - G \leq |k| \leq m} \lambda_{k, m}^{\Psi} \Psi_k^{(\gamma + G)}
      (\theta), &  & |m| > G,
    \end{array} \right.
  \end{array}
\end{equation}
i.e. $\Psi_m^{(\gamma)}$ is a linear combination of at most $2 G + 1$
functions $\Psi_k^{(\gamma + G)}$. Note that the Fourier relation
(\ref{eq:fourier-fun-connection}) parallels (\ref{eq:fourier-mode-connection})
in exactly the same way that the Jacobi relations (\ref{eq:jdemotiona}) --
(\ref{eq:jdemotionb}) parallel (\ref{eq:jpromotiona}) -- (\ref{eq:jpromotionb}).

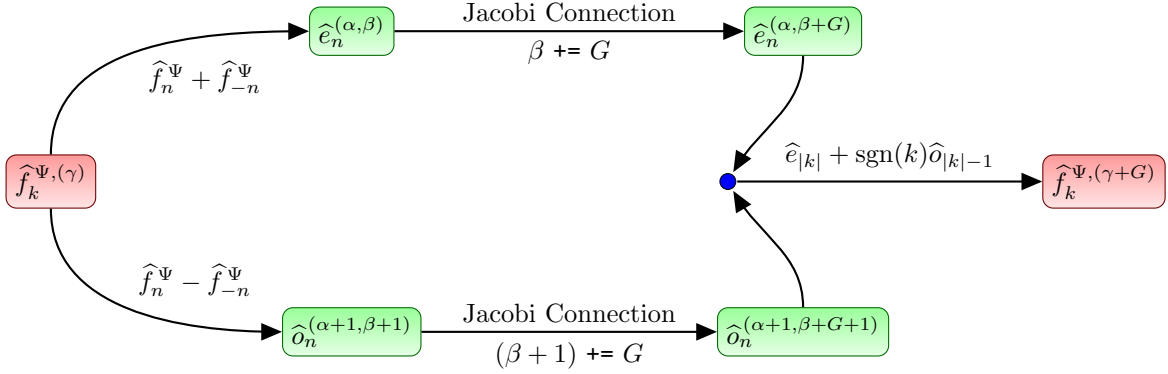
\begin{figure}[tbp]
 \begin{center}
  \input{tikz/part1/PsiPsiConnect.tex}
  \caption{Illustration of steps taken to perform $\Psi$-$\Psi$
  connections. The operator {\tt +=} is the addition-assignment operator. \label{fig:fourier-connection}}
  \end{center}
\end{figure}

We now illustrate how to calculate the Szeg$\ddot{\text{o}}$-Fourier
connection coefficients $\lambda^{\Psi}$ in Proposition
\ref{prop:Psi-connection} from the Jacobi coefficients $\lambda^P$. In the
following, we make use of the notation:
\[ \begin{array}{lllll}
     n \assign |k| - 1, &  & \alpha = - \frac{1}{2}, &  & \beta = \gamma -
     \frac{1}{2} .
   \end{array} \]
From the definition of $\Psi_k^{(\gamma)}$ in (\ref{eq:Psi-def}) we have
\[ \begin{array}{l}
     \left. \begin{array}{rcl}
       \tilde{P}_{n + 1}^{(\alpha, \beta)} & = & \Psi^{(\gamma)}_k + \Psi_{-
       k}^{(\gamma)}\\
       &  & \\
       \tilde{P}_n^{(\alpha + 1, \beta + 1)} & = & \Psi_{|k|}^{(\gamma)} -
       \Psi_{- |k|}^{(\gamma)}
     \end{array} \right\} n \geq 0,\\
     \\
     \hspace{2cm} \begin{array}{lll}
       \tilde{P}_0^{(\alpha, \beta)} & = & \sqrt{2} \Psi^{(\gamma)}_0 .
     \end{array}
   \end{array} \]
Therefore, from the modes $\hat{f}_k^{(\gamma)}$ we can derive two sets of
Jacobi modes:
\[ \begin{array}{rclll}
     \hat{e}^{(\alpha, \beta)}_n & = & \hat{f}_n^{(\gamma)} + \hat{f}_{-
     n}^{(\gamma)}, &  & n \geq 1,\\
     &  &  &  & \\
     \hat{o}^{(\alpha + 1, \beta + 1)}_n & = & \hat{f}^{(\gamma)}_{n + 1} -
     \hat{f}_{- n - 1}^{(\gamma)}, &  & n \geq 0,\\
     &  &  &  & \\
     \hat{e}^{(\alpha, \beta)}_0 & = & \sqrt{2}  \hat{f}_0^{(\gamma)} . &  & 
   \end{array} \]
The Jacobi modes $\hat{e}_n$ are modes in an expansion in polynomials
$\tilde{P}_n^{(\alpha, \beta)}$ and the modes $\hat{o}_n$ are for an expansion
in $\tilde{P}_n^{(\alpha + 1, \beta + 1)}$. With these modes in hand, we can
use the Jacobi connection coefficients to promote the coefficients using
Proposition \ref{lemma:jacobi-connection}.
\[ \begin{array}{lllllll}
     \hat{e}_n^{(\alpha, \beta + G)} & = & \sum_{m = 0}^G \lambda^P_{n, n + m}
     \hat{e}_{n + m}^{(\alpha, \beta)}, & \text{where} & \lambda^P =
     \lambda^P \left( \alpha, \beta, 0, G), \right. &  & n \geq 0\\
     &  &  &  &  &  & \\
     \hat{o}_n^{(\alpha + 1, \beta + G + 1)} & = & \sum_{m = 0}^G
     \lambda^P_{n, n + m}  \hat{o}_{n + m}^{(\alpha + 1, \beta + 1)}, &
     \text{where} & \lambda^P = \lambda^P (\alpha + 1, \beta + 1, 0, G), &  &
     n \geq 0.
   \end{array} \]
Finally we redistribute the modes back into Szeg$\ddot{\text{o}}$-Fourier form
to yield what we desired:
\[ \begin{array}{lllll}
     \hat{f}_n^{(\gamma + G)} & = & \frac{1}{2} \left[ \hat{e}_n^{(\alpha,
     \beta + G)} + \hat{o}_{n - 1}^{(\alpha + 1, \beta + 1 + G)} \right], &  &
     n \geq 1\\
     &  &  &  & \\
     \hat{f}_{- n}^{(\gamma + G)} & = & \frac{1}{2} \left[ \hat{e}_n^{(\alpha,
     \beta + G)} - \hat{o}_{n - 1}^{(\alpha + 1, \beta + 1 + G)} \right], &  &
     n \geq 1\\
     &  &  &  & \\
     \hat{f}_0^{(\gamma + G)} & = & \frac{\hat{e}_0^{(\alpha, \beta +
     G)}}{\sqrt{2}} . &  & 
   \end{array} \]
The whole procedure is illustrated graphically in Figure
\ref{fig:fourier-connection}.
We may explicitly write the connections as:
\begin{align*}
     \hat{f}_k^{(\gamma + G)} & = \sum_{m = 0}^G \frac{1}{2} \left[
     \lambda^{P ; (\alpha, \beta)}_{|k|, |k| + m} + \tmop{sgn} (k) \lambda^{P
     ; (\alpha + 1, \beta + 1)}_{|k|, |k| + m} \right] \hat{f}_{|k| +
     m}^{(\gamma)} + \\
     &  \\
     &  \hspace{1cm} \sum_{m = 0}^G \frac{1}{2} \left[ \lambda^{P ;
     (\alpha, \beta)}_{|k|, |k| + m} - \tmop{sgn} (k) \lambda^{P ; (\alpha +
     1, \beta + 1)}_{|k|, |k| + m} \right] \hat{f}_{- |k| - m}^{(\gamma)}, 
     \hspace{1cm} |k| \geq 1\\
     &  \\
     \hat{f}_0^{(\gamma + G)} & = \frac{1}{\sqrt{2}} \sum_{m = 0}^G
     \lambda^{P ; (\alpha, \beta)}_{0, m}  \hat{f}_m^{(\gamma)} +
     \frac{1}{\sqrt{2}} \sum_{m = 0}^G \lambda_{0, m}^{P ; (\alpha, \beta)}
     \hat{f}_{- m}^{(\gamma)} . 
\end{align*}
Therefore we have an explicit expression for the Szeg$\ddot{\text{o}}$-Fourier
connection coefficients in (\ref{prop:Psi-connection}):
\begin{equation}
  \label{eq:lambda-Psi} \begin{array}{lll}
    \lambda^{\Psi}_{k, \pm \left( |k| + m \right)} & = & \left\{
    \begin{array}{lll}
      \frac{1}{2} \left[ \lambda^{P ; (\alpha, \beta)}_{|k|, |k| + m} \pm
      \tmop{sgn} (k) \lambda^{P ; (\alpha + 1, \beta + 1)}_{|k|, |k| + m}
      \right], &  & |k| \geq 1\\
      &  & \\
      \frac{1}{\sqrt{2}} \lambda^{P ; (\alpha, \beta)}_{0, m}, &  & k = 0
    \end{array} \right.
  \end{array}
\end{equation}
Of course, owing to observation (\ref{eq:Psi-connection-orthogonality}), the
above equation restricts $0 \leq m \leq G$. As mentioned, this connection relation is also
valid for converting an expansion in the functions $\Phi_k^{(s)} (x)$ to one
in the functions $\Phi_k^{(s + S)} (x)$ for $S \in \mathbbm{N}$ since the
modes for these two expansions are the same. Let $f (\theta)$ be given
and define $g (x) = f (\theta (x))$. Then for all $\gamma > \frac{1}{2}$:
\[ \hat{f}^{\Psi, (\gamma)}_k \circeq \left\langle f, \Psi_k^{(\gamma)}
   \right\rangle_{w_{\theta}^{(\gamma, 0)}} \equiv \left\langle g,
   \Phi_k^{(\gamma + 1)} \right\rangle_{w_x^{(\gamma + 1, 0)}} \circeq
   \hat{g}_k^{\Phi, (\gamma + 1)} \]
This completes the $\Psi$-$\Psi$ connection problem. The
reverse connection problem (converting $\hat{f}_k^{\Psi, (\gamma + G)}$ modes
to $\hat{f}_k^{\Psi, (\gamma)}$ modes) is solved by reversing the
above procedure (all steps are invertible) and use of the fact that the forward Jacobi connection problem with
integral separation is banded upper-triangular and thus the backward connection
is $\mathcal{O} (N)$ calculable
sequentially via back-substitution. See {\cite{narayan2009}}.

We have determined how to quickly and exactly accomplish the
connection problems for the unweighted functions
\[ \begin{array}{lll}
     \sum_k \hat{f}_k^{\Psi, (\gamma)} \Psi_k^{(\gamma)} (\theta) &
     \longleftrightarrow & \sum_k \hat{f}_k^{\Psi, (\gamma + G)}
     \Psi_k^{(\gamma + G)} (\theta),\\
     &  & \\
     \sum_k \hat{g}_k^{\Phi, (s)} \Phi_k^{(s)} (x) & \longleftrightarrow &
     \sum_k \hat{g}_k^{\Phi, (s + S))} \Phi_k^{(s + S)} (x),
   \end{array} \]
in $\mathcal{O} (N)$ time where $N$ is the total number of modes when $S, G \in
\mathbbm{Z}$. These connections can be performed by utilizing the connection
coefficients in (\ref{eq:lambda-Psi}) along with the sparse connection result
of Proposition \ref{prop:Psi-connection}. For $S, G \nin \mathbbm{Z}$, there
is no sparse connection result for the modes, and so while the connection
coefficients $\lambda_{k, l}^{\Psi}$ can still be calculated based on known
connection coefficients for Jacobi polynomials, the coefficients do not
terminate finitely, and it is more expensive (that is, more costly than
$\mathcal{O}(N)$) to change $s$ or $\gamma$.

We have not described the details of how this $\Psi$-$\Psi$ connection
problem relates to the two issues presented at the beginning of this section
(i.e., using the FFT and modification of $s$ for the weighted functions
$\phi^{(s)} (x)$). The problem of using the FFT we will postpone until Part
\Rmnum{2}, which describes computational issues. In Section
\ref{sec:propjacobi-connections-smod} we will describe a method for
modification of the decay parameter $s$, for which the connection process
described in this section is an integral part.

\subsubsection{The $\Psi$-$\psi$ Connection Problem}

\label{sec:propjacobi-useless-connection}We now consider the following
problem: let $f \in L^2 \left( [- \pi, \pi], \mathbbm{C} \right)$. We assume
$\gamma \geq 0$ and consider two expansions:
\[ \begin{array}{lll}
     f (\theta) & = & \sum_{k \in \mathbbm{Z}} \hat{f}^{\Psi}_k
     \Psi_k^{(\gamma)} (\theta),\\
     &  & \\
     f (\theta) & = & \sum_{k \in \mathbbm{Z}} \hat{f}_k^{\psi}
     \psi_k^{(\gamma)} (\theta) .
   \end{array} \]
The modal coefficients are defined in the following way:
\[ \begin{array}{lll}
     \hat{f}_k^{\Psi} & = & \left\langle f, \Psi_k^{(\gamma)}
     \right\rangle_{w_{\theta}^{(\gamma)}},\\
     &  & \\
     \hat{f}_k^{\psi} & = & \left\langle f, \psi_k^{(\gamma)} \right\rangle .
   \end{array} \]
We assume that the modal coefficients for the uppercase (unweighted function)
expansion are known and that we wish to determine the lowercase modes
$\hat{f}^{\psi}$. From the definitions of the modal coefficients, it is clear
that we can rewrite the lowercase modes as
\[ \begin{array}{lll}
     \hat{f}_k^{\psi} & = & \left\langle f, \psi_k^{(\gamma)} \right\rangle\\
     &  & \\
     & = & \left\langle f \left[ \sqrt[\asterisk]{w_{\theta}^{(\gamma)}}
     \right]^{- \gamma}, \Psi_k^{(\gamma)}
     \right\rangle_{w_{\theta}^{(\gamma)}} .
   \end{array} \]
That is, the modal coefficients for the lowercase basis are identical to modal
coefficients of a different function for the uppercase basis. To see how this
helps us, we make a small digression; recall \eqref{eq:phase-shift-rewrite}
and define
\begin{equation}
  \label{eq:f-g-definition} g (\theta) \assign f \left[
  \sqrt[\asterisk]{w_{\theta}^{(\gamma)}} \right]^{- \gamma} = f \times \left[
  \frac{\sqrt{2}}{i \left( 1 + e^{- i \theta} \right)} \right]^{\gamma} .
\end{equation}
Suppose that $\gamma = G \in \mathbbm{N}_0$ and that we can somehow find the
modal coefficients
\[ \hat{g}^{\Psi, (0)}_k = \left\langle g, \Psi_k^{(0)} \right\rangle . \]
Then we can use the $\Psi$-$\Psi$ connection problem outlined in
Section \ref{sec:propjacobi-connections-fourier} to accurately and efficiently
determine the modal coefficients $\hat{g}_k^{\Psi}$ for $\gamma = G$ due to
the sparse connection. To see how we can find the modal coefficients
$\hat{g}_k^{\Psi, (0)}$, assume that we have the modal coefficients
$\hat{f}_k^{\Psi, (0)}$. Then (\ref{eq:f-g-definition}) implies that
\begin{equation}
  \label{eq:fhat-ghat-relation} \sum_{m = 0}^G \left( \begin{array}{l}
    G\\
    m
  \end{array} \right)  \hat{g}_{k + m}^{\Psi, (0)} = \hat{f}_k^{\Psi, (0)}
  \left( \frac{\sqrt{2}}{i} \right)^G .
\end{equation}
If we assume a finite expansion so that $\hat{g}_k = 0$ for $|k| > 2 N + 1$,
then we can solve (\ref{eq:fhat-ghat-relation}) via back-substitution.
Note that determining each coefficient costs $\mathcal{O} (G)$
operations, independent of $N$; this is a similar operation count to the
$\Psi$-$\Psi$ connection cost.

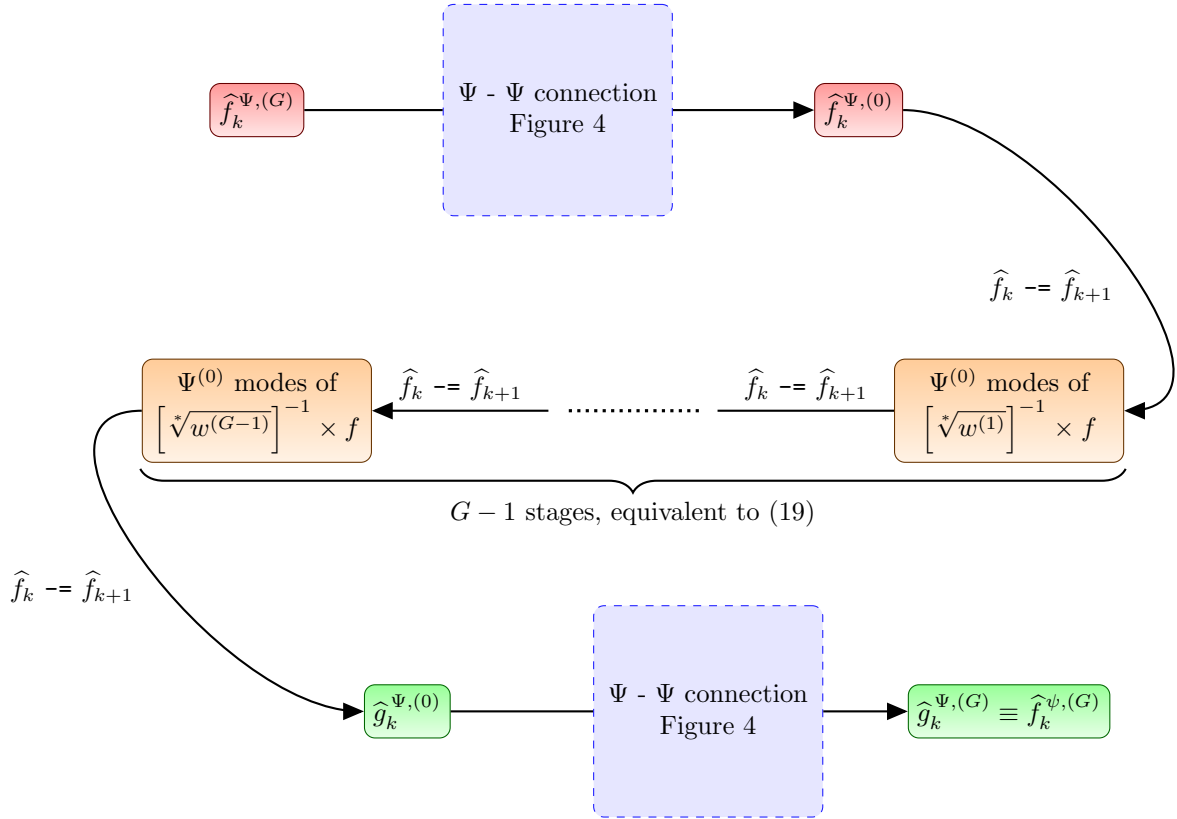
\begin{figure}[tbp]
  \begin{center} 
  \input{tikz/part1/PsipsiConnect}
  \caption{Flowchart representation of a $\Psi$-$\psi$
  connection. The operator {\tt -=} is the subtraction-assignment operator. \label{fig:Psi-psi-connection}}
  \end{center}
\end{figure}

Finally, we must obtain $\hat{f}_k^{\Psi, (0)}$ from the given input
$\hat{f}_k^{\Psi, (G)}$. However, this is another $\Psi$-$\Psi$ connection
(albeit in reverse). Therefore, the three steps to take us from
$\hat{f}_k^{\Psi, (G)}$ modes to $\hat{f}_k^{\psi, (G)}$ modes are
\begin{enumerate}
  \item Compute $\hat{f}_k^{\Psi, (0)}$ from $\hat{f}_k^{\Psi, (G)}$, which is
  a (backward) $\Psi$-$\Psi$ connection
  
  \item Compute $\hat{g}_k^{\Psi, (0)}$ from $\hat{f}_k^{\Psi, (0)}$ using
  (\ref{eq:fhat-ghat-relation}).
  
  \item Compute $\hat{f}_k^{\psi, (G)} \equiv \hat{g}_k^{\Psi, (G)}$ from
  $\hat{g}_k^{\Psi, (0)}$, a (forward) $\Psi$-$\Psi$ connection.
\end{enumerate}
This is illustrated in Figure \ref{fig:Psi-psi-connection}. For an expansion
with $N$ modes, all three steps have $\mathcal{O} (N \nonesep G)$ cost
asymptotically. The backward connection problem (determining $\hat{f}^{\Psi,
(G)}$ from $\hat{f}^{\psi, (G)}$) is also computable in $\mathcal{O} \left( NG
\right)$ operations, and is accomplished by reversing the above operations.

Note that if $\gamma \nin \mathbbm{N}_0$ then all of these steps break down:
the $\Psi$-$\Psi$ connection is not sparse, and (\ref{eq:fhat-ghat-relation}) is
not valid since $\gamma$ is not an integer in (\ref{eq:f-g-definition}).

Performing these modal connections on the real line for expansions in
$\Phi^{(s)} (x)$ and $\phi^{(s)} (x)$ is equivalent, except one must assign
$\gamma \assign s - 1$ and then proceed as outlined above.

This particular connection problem is not necessarily useful explicitly since
in many of our applications, we will have direct access to $\hat{f}^{\Psi,
(0)}$, but each of the pieces necessary for this computation are used
extensively both in modification of the decay parameter $s$ and application of
the FFT.

\subsubsection{Modification Of $s$: The $\psi$-$\psi$ Connection}

\label{sec:propjacobi-connections-smod}
\begin{figure}[tbp]
  \begin{center}
  \input{tikz/part1/psipsiConnect}
  \caption{Flowchart of operations for modification of $s$. The operator {\tt
  -=} is the subtraction-assignment operator.
  \label{fig:s-modification}}
  \end{center}
\end{figure}
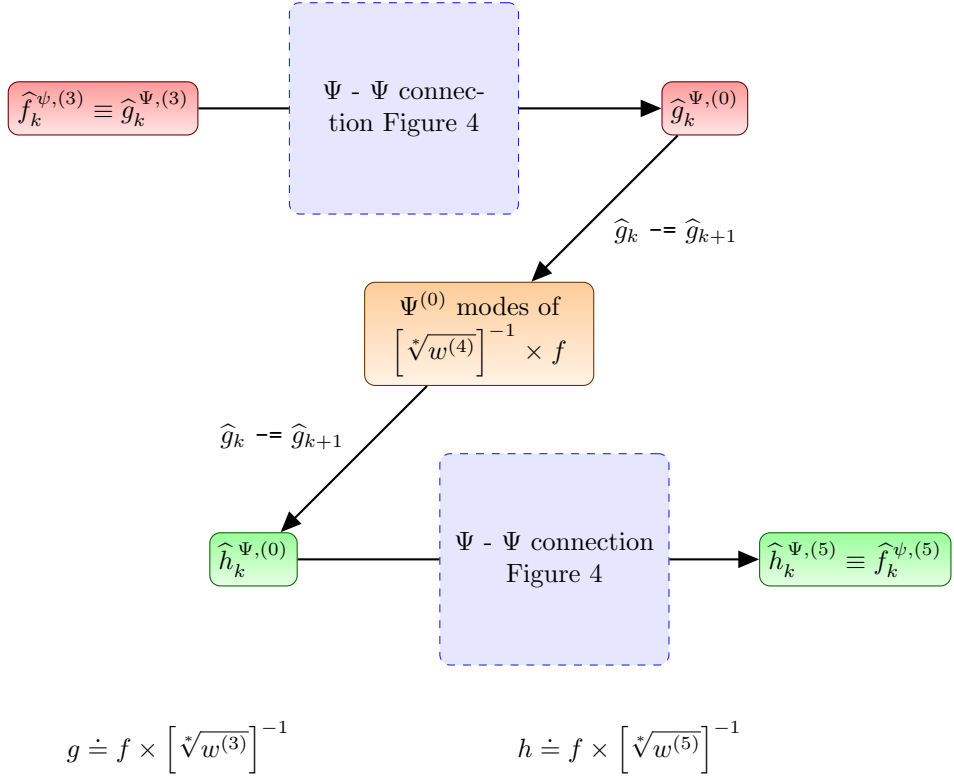

We have now developed the necessary
tools for the modification of $s$, i.e., the $\psi$-$\psi$
connection problem. We assume that $G, F \in \mathbbm{N}$ and that we know
connection coefficients of some function $f \in L^2$ for an expansion in
$\psi^{(F)}$, and wish to obtain the coefficients for a $\psi^{(G)}$
expansion. The whole procedure can be accomplished in three steps:
\begin{enumerate}
  \item Obtain expansion coefficients for $f \times \left[
  \sqrt[\asterisk]{w^{(F)}_{\theta}} \right]^{- 1}$ in the $\Psi^{(0)}$
  
  ($\Psi$-$\psi$ connection)
  
  \item Obtain expansion coefficients for $f \times \left[
  \sqrt[\asterisk]{w_{\theta}^{(G)}} \right]^{- 1}$ in the $\Psi^{(0)}$
  
  (Fourier connection)
  
  \item Obtain the sought expansion coefficients of $f$ in the $\psi^{(G)}$
  
  ($\Psi$-$\psi$ connection)
\end{enumerate}
Step 2 is easily performed using a version of (\ref{eq:fhat-ghat-relation}) by
noting the relation between $f \times \left[
\sqrt[\asterisk]{w^{(G)}_{\theta}} \right]^{- 1}$ and $f \times \left[
\sqrt[\asterisk]{w_{\theta}^{(F)}} \right]^{- 1}$ with knowledge of the
canonical Fourier expansion coefficients ($\Psi^{(0)} (\theta)$). This is
shown in Figure \ref{fig:s-modification} for the special case $F=3$, $G=5$. 

Note that this particular connection problem is very amenable to an
FFT+collocation approach whereas the algorithm we have laid out is a
`Galerkin' approach. The problem with the collocation approach is that it
requires $\mathcal{O} \left( N \log N) \right.$ operations with two FFT's,
whereas the above algorithm requires only $\mathcal{O} \left( NG) \right.$
steps.

As with the $\Psi$-$\psi$ connection of the previous section, if either the
starting parameter $F$ or the target parameter $G$ are not integers, then
this procedure cannot be used: the core of the fast algorithm is the ability to
obtain the canonical Fourier modes, which cannot be done efficiently if the
decay parameters are not integers. 

\subsection{Quadrature}

\label{sec:wiener-quadrature}We now turn to quadrature rules that will compute
integrals over the real line. We adopt the following notation: the pair
$\left\{ r^{(\alpha, \beta)}_n, \omega^{(\alpha, \beta)}_n \right\}_{n = 1}^N$
denotes the $N$-point Gauss-quadrature for the Jacobi polynomial of class
$(\alpha, \beta)$, i.e.,
\[ \int_{- 1}^1 f (r) w_r^{(\alpha, \beta)} \mathd r = \sum_{n = 1}^N f \left(
   r_n^{(\alpha, \beta)} \right) \omega_n^{(\alpha, \beta)}, \hspace{1cm}
   \forall f \in \mathcal{B}_{2 N-1} , \]
where $\mathcal{B}_{2 N-1}$ is the space of polynomials of degree $2 N -1$ or
less. We suppress the dependence of $r_n^{(\alpha, \beta)}$ and $\omega_n^{(\alpha,
\beta)}$ on $N$. We also denote $\left\{ r_n^{(\alpha, \beta) ; \tmop{GR}},
\omega_n^{(\alpha, \beta) ; \tmop{GR}} \right\}_{n = 1}^N$ as the $N$-point
Gauss-Radau quadrature with the fixed node $r_N^{(\alpha, \beta) ; \tmop{GR}}
\equiv 1$. We assume for clarity of presentation
that the nodes are ordered by $n$, e.g. $r_{n - 1}^{(\alpha, \beta)} <
r_n^{(\alpha, \beta)}$.

With the goal that we wish to develop quadrature rules for the
infinite line, we will take pains to develop quadrature rules in
$\theta$-space that do not have nodes at $\theta = \pm \pi$, which map to $x =
\pm \infty$. We use the Jacobi-Gauss quadrature rules as the building blocks
for our generalized Fourier quadrature rules.

Suppose we wish to construct an $N$-point quadrature rule associated with the
functions $\Psi_k^{(\gamma)} (\theta)$. If $N$ is even, then define
\begin{equation}
  \label{eq:quad-theta-even} \theta^{(\gamma)}_n = \left\{ \begin{array}{lll}
    - \arccos \left( r_n^{(- 1 / 2, \gamma - 1 / 2)} \right), &  & 1 \leq n
    \leq \frac{N}{2}\\
    &  & \\
    - \theta_{N + 1 - n}^{(\gamma)}, &  & \frac{N}{2} + 1 \leq n \leq N,
  \end{array} \right.
\end{equation}
where $r_n^{(\alpha, \beta)}$ comes from an $\frac{N}{2}$-point quadrature
rule, and
\[ \Omega_n^{(\gamma)} = \left\{ \begin{array}{lll}
     \omega_n^{(- 1 / 2, \gamma - 1 / 2)}, &  & 1 \leq n \leq \frac{N}{2}\\
     &  & \\
     \Omega_{N + 1 - n}^{(\gamma)}, &  & \frac{N}{2} + 1 \leq n \leq N,
   \end{array} \right. \]
and $\omega_n^{(\alpha, \beta)}$ comes from an $\frac{N}{2}$ -point quadrature
rule.

If $N$ is odd, then define
\begin{equation}
  \label{eq:quad-theta-odd} \theta^{(\gamma)}_n = \left\{ \begin{array}{lll}
    - \arccos \left( r_n^{(- 1 / 2, \gamma - 1 / 2) ; \tmop{GR}} \right), &  &
    1 \leq n \leq \frac{N + 1}{2}\\
    &  & \\
    \theta_{N + 1 - n}^{(\gamma)}, &  & \frac{N + 3}{2} \leq n \leq N,
  \end{array} \right.
\end{equation}
where $r_n^{(\alpha, \beta) ; \tmop{GR}}$ comes from an $\frac{N +
1}{2}$-point quadrature rule, and
\[ \Omega_n^{(\gamma)} = \left\{ \begin{array}{lll}
     \omega_n^{(- 1 / 2, \gamma - 1 / 2) ; \tmop{GR}}, &  & 1 \leq n \leq
     \frac{N - 1}{2}\\
     &  & \\
     2 \omega_n^{(- 1 / 2, \gamma - 1 / 2) ; \tmop{GR}}, &  & n = \frac{N +
     1}{2}\\
     &  & \\
     \Omega_{N + 1 - n}^{(\gamma)}, &  & \frac{N + 3}{2} \leq n \leq N.
   \end{array} \right. \]

\begin{figure}[tbp]
  \begin{center}
  \input{tikz/part1/GaussConstruct}
  \caption{Construction of Gauss-type quadrature for generalized Fourier
  functions. The new quadrature rules are symmetric combinations of
  Jacobi-Gauss-type quadrature rules. The constructions shown are accurate
  node locations for $\gamma = 5$.}
  \label{fig:fourier}
  \end{center}
\end{figure}
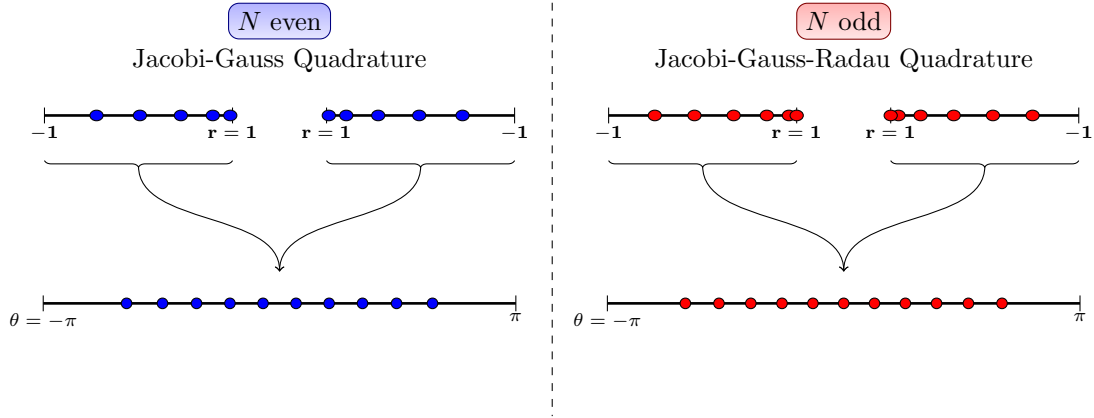

\noindent For graphical descriptions of the above formulae, see Figure
\ref{fig:fourier}.  We have used Jacobi-Gauss rules for $N$ even and
Jacobi-Gauss-Radau rules for $N$ odd. By construction, when $N$ is odd,
$\theta_{\frac{N + 1}{2}}^{(\gamma)} = 0$ due to the Gauss-Radau rule
requirement that $r_{\frac{N + 1}{2}}^{(\alpha, \beta) ; GR} = 1$. The quadrature
rules derived above have no nodes at $\theta = \pm \pi$ (since there are no
Jacobi-Gauss, or Jacobi-Gauss-Radau nodes at $r = - 1$) and are symmetric rules
for any $\gamma$. Thus they are always exact for any odd function. It is not
difficult to show the following result:

\begin{proposition}
  For $N$ even, the $N$-point quadrature rule $\left\{ \theta_n^{(\gamma)},
  \Omega_n^{(\gamma)} \right\}_{n = 1}^N$ satisfies
  \[ \int_{- \pi}^{\pi} e^{i \nonesep k \nonesep \theta} w_{\theta}^{(\gamma,
     0)} \mathd \theta = \sum_{n = 1}^N e^{i \nonesep k \nonesep
     \theta_n^{(\gamma)}} \left( \theta_n^{(\gamma)} \right)
     \Omega_n^{(\gamma)}, \hspace{1cm} |k| \leq N - 1. \]
  When $N$ is odd, the quadrature rule satisfies
  \[ \begin{array}{l}
       \int_{- \pi}^{\pi} e^{i \nonesep k \nonesep \theta}
       w_{\theta}^{(\gamma, 0)} \mathd \theta = \sum_{n = 1}^N e^{i \nonesep k
       \nonesep \theta_n^{(\gamma)}} \left( \theta_n^{(\gamma)} \right)
       \Omega_n^{(\gamma)}, \hspace{1cm} |k| \leq N.
     \end{array} \]
\end{proposition}

\noindent The degeneracy in the quadrature rule for $N$ even is exactly of the same
nature as the degeneracy in the canonical equispaced Fourier quadrature rule
for an even number of grid points \cite{hesthaven2007}. If $\gamma = 0$ the rule $\left\{
\theta_n^{(0)}, \Omega_n^{(0)} \right\}_{n = 1}^N$ is exactly the same as the
equispaced Fourier quadrature rule, symmetric about $\theta = 0$. The
quadrature rule $\left\{ \theta_n^{(0)}, \Omega_n^{(0)} \right\}_{n = 1}^N$
can be used to integrate against the weight function $w_{\theta}^{(\gamma,
0)}$ when $\gamma \in \mathbbm{N}$ since in this case the weight is itself a
trigonometric polynomial.

In order to determine a quadrature rule to integrate the weighted functions
$\psi_k^{(\gamma)} (\theta)$, we can augment the weights $\Omega_n^{(\gamma)}$
to contain information about the weight function. This can be summed up in the
following result:

\begin{corollary}
  The even $N$-point quadrature rule $\left\{ \theta_n^{(\gamma)},
  \omega_n^{(\gamma)} \right\}_{n = 1}^N$, where $\omega_n^{(\gamma)} \assign
  w_{\theta}^{(- \gamma, 0)} \left( \theta_n^{(\gamma)} \right)
  \Omega_n^{(\gamma)}$ satisfies
  \[ \int_{- \pi}^{\pi} \psi_k^{(\gamma)} \overline{\psi_l^{(\gamma)}} \mathd
     \theta = \sum_{n = 1}^N \psi_k^{(\gamma)} \left( \theta_n^{(\gamma)}
     \right) \overline{\psi_l^{(\gamma)}} \left( \theta_n^{(\gamma)} \right)
     \omega_n^{(\gamma)}, \hspace{1cm} |k| + |l| \leq N - 1 \]
\end{corollary}

Multiplying $\Omega_n^{(\gamma)}$ by the inverse of the weight $w_{\theta}^{(-
\gamma, 0)}$ is mathematically not a problem since none of the
$\theta_n^{(\gamma)}$ are equal to $\pm \pi$, where the weight $w_{\theta}^{(-
\gamma, 0)}$ is singular. Note that since the functions $\Phi_k^{(s)} (x)$ are
just a mapping of the functions $\Psi_k^{(\gamma)} (\theta)$, the
quadrature rule $\left\{ x \left( \theta_n^{(s - 1)} \right), \Omega_n^{(s -
1)} \right\}_{n = 1}^N$, which has nodal values over $\mathbbm{R}$, can be
used to integrate the functions $\Phi_k^{(s)} (x)$ over the real line.
Similarly, the rule $\left\{ x \left( \theta_n^{(s - 1)} \right), \omega_n^{(s
- 1)} \right\}_{n = 1}^N$ can be used to integrate Galerkin products of the
generalized Wiener functions $\phi_k^{(s)} (x)$ over the real line.

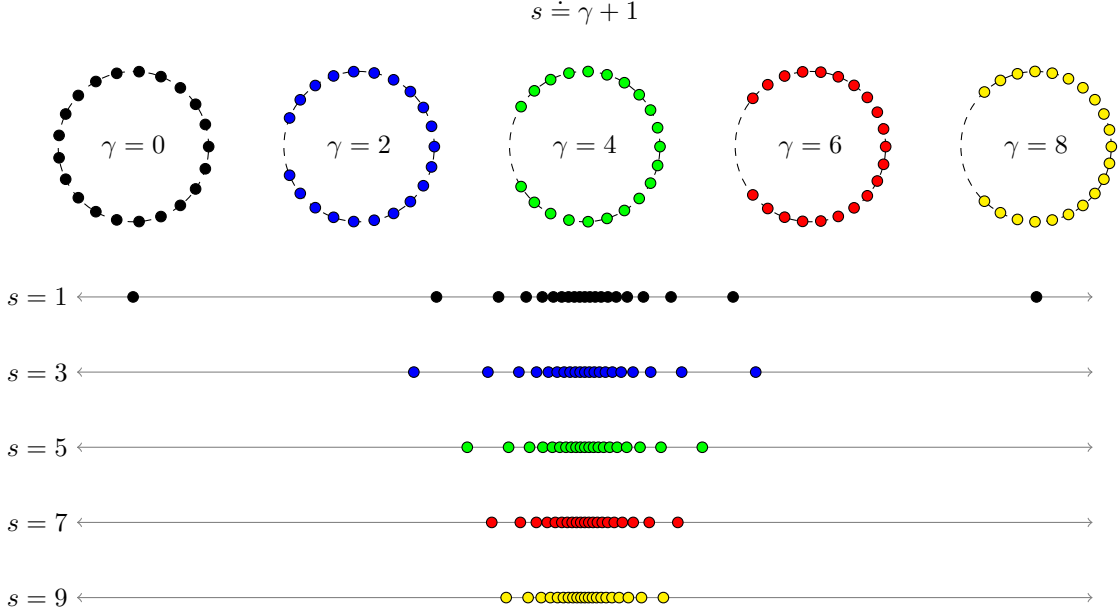
\begin{figure}[tbp]
  \begin{center}
  \input{tikz/part1/quadpoints}
  \caption{(Top) Plots of the Fourier quadrature nodes on the unit circle
  generated with equation (\ref{eq:quad-theta-odd}), $N = 21$. (Bottom) The
  resulting quadrature nodes on the real line. The scale on the real line is
  $|x| \leq 15$.}
  \label{fig:quad-point-plots}
  \end{center}
\end{figure}

For various $\gamma / s$ we graphically depict the location of the quadrature
nodes for $N = 21$ in Figure \ref{fig:quad-point-plots} on the unit circle $z
\in \mathbbm{T}$ and on the real line. Note that as we increase $\gamma$ the quadrature nodes
become more and more concentrated towards $z = 1$ ($\theta = 0$). On the real
line, this manifests itself as higher concentration near $x = 0$ which,
although rectifiable via an affine mapping, is suboptimal if one wishes to
resolve functions away from $x = 0$. Note that the tendency of Jacobi-Gauss
nodes to become more equidistant on $[- 1, 1]$ as $\beta$ (i.e. $\gamma$ or
$s$) is increased {\cite{hesthaven2007}} also suggests that these generalized
quadrature rules for large $\gamma$ or $s$ will not be as good as the the ones
for smaller $\gamma$ or $s$ since equidistant nodes are bad for
finite-interval polynomial interpolation. In addition, when $\gamma = 0$, we
can use these (equidistant) quadrature nodes to employ the FFT for modal-nodal
transformations.

\subsection{The Stiffness Matrix}

\label{sec:propjacobi-stiffness}
In many applications to differential equations it is necessary to express the
derivative of a basis function as a linear combination of basis functions. We
devote this section to this endeavor. We define entries of the stiffness
matrix as
\[ S^{\phi}_{k, l} = \left\langle \phi_k^{(s)}, \frac{\mathd}{\mathd x}
   \phi_l^{(s)} \right\rangle . \]
For the generalized Wiener rational functions, the following result can be
proven:

\begin{theorem}
  \label{thm:stiffness-brief}Let $S^{\phi}$ denote the $N \times N$ stiffness
  matrix for the weighted Wiener rational functions $\phi_k^{(s)}$. $S^{\phi}$
  satisfies the following properties for any $s > \frac{1}{2}$:
  \begin{enumerate}
    \item $S^{\phi}$ is skew-Hermitian, i.e. $S^{\phi}_{k, l} = -
    \overline{S^{\phi}_{l, k}}$
    
    \item $S^{\phi}$ is sparse with entries only on the super-, sub-, and
    main sinister and dexter diagonals: define
    \[ \begin{array}{lllll}
         k^{\vee} \assign \tmop{sgn} (k) \left( |k| - 1 \right) = k -
         \tmop{sgn} (k), &  &  &  & \left. k^{\wedge} \assign \tmop{sgn} (k)
         (|k| + 1 \right) = k + \tmop{sgn} (k) .
       \end{array} \]
    Then
    \[ \frac{\mathd \phi_k^{(s)} (x)}{\mathd x} = \sum_{l \in \left\{ \pm
       k^{\vee}, \pm k, \pm k^{\wedge} \right\}} \tau_{k, l}^{(s)}
       \phi_l^{(s)} (x), \]
    for some purely imaginary constants $\tau_{k, l}^{(s)}$. In
    other words,
    \[ S^{\phi}_{k, l} = 0, \hspace{1cm} l \nin \left\{ \pm k^{\vee}, \pm k,
       \pm k^{\wedge} \right\} . \]
    \item The spectral radius of $S^{\phi}$ satisfies
    \[ \rho(S^{\phi}) \leq N + 5 s. \]
  \end{enumerate}
\end{theorem}

\noindent The proof of Theorem \ref{thm:stiffness-brief} is quite tedious, so we only
sketch the main points. Details are given in Appendix \ref{app:stiff}.

\begin{proof}
  Property 1 can easily be deduced by using integration by parts and noting
  that the functions $\phi_k^{(s)} (x)$ decay to zero as $|x| \rightarrow
  \infty$.

  Property 2 is a highly nontrivial result that is provable using several
  properties of Jacobi Polynomials. We refer the reader to
  {\cite{narayana:2009:1}}. Most of the calculations are straightforward once
  a list of Jacobi Polynomial properties has been compiled. However, there are
  some difficulties whose resolutions rely on a couple of fortuitous
  properties: first, that $\frac{\mathd x (\theta)}{\mathd \theta} (\theta) =
  1 + \cos \theta$, i.e. that the map we have chosen to take $\theta
  \rightarrow x$ has a Jacobian with a particular form. Second, that
  \[ \frac{\mathd}{\mathd \theta} \left[ \left( \sin \theta \right)
     \tilde{P}_n^{(\alpha + 1, \beta + 1)} \left( \cos \theta \right) \right]
  \]
  is a sparse combination of $\tilde{P}_n^{(\alpha, \beta + 1)} \left( \cos
  \theta \right)$, which is not a trivial result; we show this by using
  brute-force calculation with the compiled list of Jacobi
  Polynomial properties.

  Property 3 can be derived from the second property. The key ingredient is
  Gerschgorin's Theorem. Using the explicit entries for the constants
  $\tau_{k, l}^{(s)}$ given in Theorem \ref{thm:stiff-entries} of Appendix
  \ref{app:stiff} we can show that for each $k$ satisfying $|k| \geq 2$ the
  following crude bounds hold:
  
  \[ \begin{array}{rcl}
       | \tau_{k, k} | & \leq & n + 2 s,\\
       &  & \\
       | \tau_{k, - k} | + | \tau_{k, k^{\vee}} | + | \tau_{k, - k^{\vee}} | +
       | \tau_{k, k^{\wedge}} | + | \tau_{k, - k^{\wedge}} | & \leq & n + 3 s
       + 2,
     \end{array} \]

  where $n \assign |k| - 1$. Gerschgorin's Theorem can now be used to define a
  region in the complex plane in which all the eigenvalues lie. By the above
  properties, this region has distance from the origin at most $2 n + 5 s +
  2$. Once we consider the necessary relationship between $n$, $k$, and $N$,
  the result is proven. (It is interesting, but not necessary, to note that
  the eigenvalues all lie on the imaginary axis due to the skew-Hermitian
  property of $S$.)
\end{proof}

\begin{remark}
  While the $\mathcal{O} (N)$ maximum eigenvalue does depend on $s$, the proportionality
  factor is empirically around 2, not $5$ as given in the theorem. See Table
  \ref{tab:stiff-eig}.
\end{remark}

\begin{figure}[tbp]
  \begin{center}
  \includegraphics[width=1.0\textwidth]{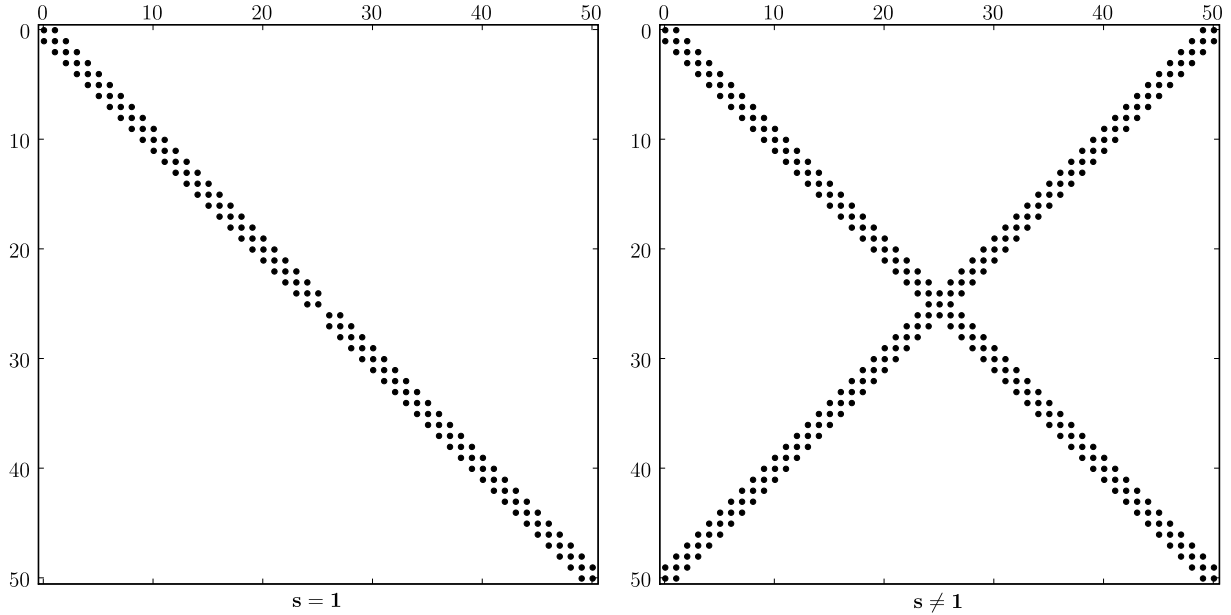}
  \caption{\label{fig:dpsi-sparse}Sparsity plots for stiffness matrices of the
  weighted Wiener rational functions $\phi_k^{(s)}$. The sparsity patterns are
  representative of property 2 in Theorem \ref{thm:stiffness-brief} for $s =
  1$ (left) and all $s \neq 1$ (right). The $s = 1$ sparsity pattern has been
  derived previously {\cite{christov1982}}, and the expressions for the
  $\tau_{k, l}$ in Appendix \ref{app:stiff} with $s = 1$ reduce to the pattern
  above.}
  \end{center}
\end{figure}

The sparsity pattern we have derived for the derivatives of these functions
(property 2 of the above theorem) is illustrated in Figure
\ref{fig:dpsi-sparse}. Note that the unweighted functions $\Phi^{(s)}(x)$ also have
a similar sparsity result; see Lemma \ref{lemma:stiff2}. However, the Fourier
functions $\Psi^{(\gamma)}(\theta)$ and $\psi^{(\gamma)}(\theta)$ do not have sparse
stiffness matrices (unless $\gamma=0$).  In addition, numerical values for the maximum
eigenvalues of the stiffness matrix (property 3) are given in Table
\ref{tab:stiff-eig}. The sparsity of the stiffness matrix is important for
fast computations of derivatives for spectral methods for solving PDEs, and
the $\mathcal{O} (N)$ maximum eigenvalue of the stiffness matrix indicates that we can
take a relatively large timesteps for time-dependent problems. Finally, the
skew-symmetry of the stiffness matrix easily leads to energy conservation for
the Galerkin discretization of hyperbolic conservation laws.

\begin{table}[h]
\setlength{\doublerulesep}{0pt}
\renewcommand{\arraystretch}{2.0}
\centering 
  \begin{tabular}{|r!{\vrule width 1pt}c|c|c|c|c|}\hline
    $s\backslash$N &  11 & 50 & 101 & 250 & 501\\\hline\hline
    0.6 & 7.31 & 43.76 & 91.50 & 237.60 & 483.75\\\hline
    1.0 & 7.99 & 44.51 & 92.28 & 238.39 & 484.54\\\hline
    6.0 & 15.96 & 53.75 & 101.81 & 248.14 & 494.40\\\hline
    $\pi^2 \approx 9.87$ & 21.72 & 60.67 & 109.05 & 255.63 & 501.99\\\hline
    15.5 & 29.73 & 70.45 & 119.40 & 266.44 & 512.99\\\hline
  \end{tabular}
  \caption{\label{tab:stiff-eig}Maximum eigenvalue of the $N \times N$
  stiffness matrix $S^\phi$ for the Wiener rational functions $\phi_k^{(s)}$. The
  results adhere to the asymptotic bound given in property 3 of Theorem
  \ref{thm:stiffness-brief}. }
\end{table}

%% file: tikz/part1/PsiPsiConnect.tex
\begin{tikzpicture}
[Psistyle/.style={shape=rectangle,rounded corners, 
                  text centered,
                  top color=\PsiModeColorTop,
                  bottom color=\PsiModeColorBottom,
                  draw=\PsiModeColorBorder},
Jacstyle/.style={shape=rectangle,rounded corners, 
                  text centered,
                  top color=\JacModeColorTop,
                  bottom color=\JacModeColorBottom,
                  draw=\JacModeColorBorder},
sumstyle/.style={fill=blue},
connectarrow/.style={->,thick,>=triangle 45},
auto]

\def\PsiModeColorTop{red!40!white}
\def\PsiModeColorBottom{red!10!white}
\def\PsiModeColorBorder{red!40!black}
\def\JacModeColorTop{green!40!white}
\def\JacModeColorBottom{green!10!white}
\def\JacModeColorBorder{green!40!black}

\def\PsiLeft{0}   
\def\JacLeft{4}   
\def\JacRight{10} 
\def\PsiRight{14} 

\def\jacsumXoff{-1} 

\def\PsiY{-1.5}   
\def\JacYSpace{2} 

\begin{scope}[]
  \node (Psi1) at (\PsiLeft,\PsiY) [Psistyle] {$\widehat{f}_k^{\,\Psi,(\gamma)}$};

  \node (jace1) at (\JacLeft,\PsiY+\JacYSpace) [Jacstyle] {$\widehat{e}_n^{\,(\alpha,\beta)}$};
  \node (jaco1) at (\JacLeft,\PsiY-\JacYSpace) [Jacstyle] {$\widehat{o}_n^{\,(\alpha+1,\beta+1)}$};

  \node (jace2) at (\JacRight,\PsiY+\JacYSpace) [Jacstyle] {$\widehat{e}_n^{\,(\alpha,\beta+G)}$};
  \node (jaco2) at (\JacRight,\PsiY-\JacYSpace) [Jacstyle] {$\widehat{o}_n^{\,(\alpha+1,\beta+G+1)}$};

  \draw[sumstyle] (\JacRight+\jacsumXoff,\PsiY) circle (3pt) node (jacsum) {};
  
  \node (Psi2) at (\PsiRight,\PsiY) [Psistyle] {$\widehat{f}_k^{\,\Psi,(\gamma+G)}$};

  \draw [connectarrow] (Psi1) to [out=90,in=180] node [swap]
        {$\widehat{f}_n^{\,\Psi} + \widehat{f}_{-n}^{\,\Psi}$} (jace1);

  \draw [connectarrow] (Psi1) to [out=-90,in=180] node 
        {$\widehat{f}_n^{\,\Psi} - \widehat{f}_{-n}^{\,\Psi}$} (jaco1);

  \draw [connectarrow] (jace1) to [out=0,in=180] node 
        {Jacobi Connection} node [swap] {$\beta$ {\tt +=} $G$} (jace2);

  \draw [connectarrow] (jaco1) to [out=0,in=180] node
        {Jacobi Connection} node [swap] {$(\beta+1)$ {\tt +=} $G$} (jaco2);

  \draw [connectarrow] (jace2) to [out=270,in=60] (jacsum);
  \draw [connectarrow] (jaco2) to [out=90,in=-60] (jacsum);
  
  \draw [connectarrow] (jacsum) to [out=0,in=180] node
        {$\widehat{e}_{|k|} + \mathrm{sgn}(k) \widehat{o}_{|k|-1}$} (Psi2);
\end{scope}
\end{tikzpicture}

%% file: tikz/part1/PsipsiConnect.tex
\begin{tikzpicture}
[Psistyle/.style={shape=rectangle,rounded corners, 
                  text centered,
                  top color=\PsiModeColorTop,
                  bottom color=\PsiModeColorBottom,
                  draw=\PsiModeColorBorder},
psistyle/.style={shape=rectangle,rounded corners, 
                  text centered,
                  top color=\psiModeColorTop,
                  bottom color=\psiModeColorBottom,
                  draw=\psiModeColorBorder},
connectstyle/.style={shape=rectangle,rounded corners, 
                  text centered,
                  top color=\ConModeColorTop,
                  bottom color=\ConModeColorBottom,
                  text width=80pt,
                  draw=\ConModeColorBorder},
figstyle/.style={shape=rectangle,rounded corners, dashed,
                 text centered,
                 top color=\figrefColorTop,
                 bottom color=\figrefColorBottom,
                 minimum height = 80,
                 minimum width = 80,
                 text width=80pt,
                 draw=\figrefColorBorder},
sumstyle/.style={fill=blue},
connectarrow/.style={->,thick,>=triangle 45},
connectline/.style={-,thick},
connectdots/.style={-,very thick,dotted},
auto]

\def\PsiModeColorTop{red!40!white}
\def\PsiModeColorBottom{red!10!white}
\def\PsiModeColorBorder{red!40!black}
\def\ConModeColorTop{orange!40!white}
\def\ConModeColorBottom{orange!10!white}
\def\ConModeColorBorder{orange!40!black}
\def\psiModeColorTop{green!40!white}
\def\psiModeColorBottom{green!10!white}
\def\psiModeColorBorder{green!40!black}
\def\figrefColorTop{blue!10!white}
\def\figrefColorBottom{blue!10!white}
\def\figrefColorBorder{blue!90!white}

\def\braceOff{1.55}  

\def\PsiX{0}   
\def\psiX{10} 
\def\PsiY{0}   
\def\psiY{-8}   

\def\connectY{-4}  
\def\connectXStart{10}  
\def\connectXEnd{0}  
\def\connectBeforeDots{6} 
\def\connectAfterDots{4} 

\def\PsiPsiFXOffset{4}  

\begin{scope}[]
  \node (FPsi1) at (\PsiX,\PsiY) [Psistyle] {$\widehat{f}_k^{\,\Psi,(G)}$};

  \node (PsiPsiBackward) at (\PsiX+\PsiPsiFXOffset,\PsiY) [figstyle] 
        {$\Psi$ - $\Psi$ connection \\Figure \ref{fig:fourier-connection}};

  \node (FPsi2) at (\PsiX+2*\PsiPsiFXOffset,\PsiY) [Psistyle] 
        {$\widehat{f}_k^{\,\Psi,(0)}$};


  \node (connectStart) at (\connectXStart,\connectY) [connectstyle] 
        {$\Psi^{(0)}_{}$ modes of $\left[\sqrt[\ast]{w^{(1)}}\right]^{-1}\times f$};

  \node (connectDotStart) at (\connectBeforeDots,\connectY) {};
  \node (connectDotEnd) at (\connectAfterDots,\connectY) {};

  \node (connectEnd) at (\connectXEnd,\connectY) [connectstyle]
        {$\Psi^{(0)}_{}$ modes of $\left[\sqrt[\ast]{w^{(G-1)}}\right]^{-1}\times f$};


  \node (GPsi1) at (\psiX-2*\PsiPsiFXOffset,\psiY) [psistyle] 
        {$\widehat{g}_k^{\,\Psi, (0)}$};

  \node (PsiPsiForward) at (\psiX-\PsiPsiFXOffset,\psiY) [figstyle] 
        {$\Psi$ - $\Psi$ connection\\ Figure \ref{fig:fourier-connection}};

  \node (GPsi2) at (\psiX,\psiY) [psistyle] 
        {$\widehat{g}_k^{\,\Psi,(G)} \equiv \widehat{f}_k^{\,\psi,(G)}$};

  \draw [connectline] (FPsi1) to [out=0,in=180] (PsiPsiBackward);
  \draw [connectarrow] (PsiPsiBackward) to [out=0,in=180] (FPsi2);

  \draw [connectarrow] (FPsi2) to [out=0,in=0] node [swap]
        {$\widehat{f}_{k}$ {\tt -=} $\widehat{f}_{k+1}$} (connectStart);

  \draw [connectline] (connectStart) to [out=180,in=0] node [swap]
        {$\widehat{f}_k$ {\tt -=} $\widehat{f}_{k+1}$} (connectDotStart);
  \draw [connectdots] (connectDotStart) to (connectDotEnd);
  \draw [connectarrow] (connectDotEnd) to [out=180,in=0] node [swap]
        {$\widehat{f}_k$ {\tt -=} $\widehat{f}_{k+1}$} (connectEnd);

  \draw[decorate,thick,decoration={brace,amplitude=9,raise=22pt}] (\connectXStart+\braceOff,\connectY) to 
       node[yshift=-30pt] {$G-1$ stages, equivalent to
       (\ref{eq:fhat-ghat-relation})} (\connectXEnd-\braceOff,\connectY);

  \draw [connectarrow] (connectEnd) to [out=180,in=180] node [swap]
        {$\widehat{f}_{k}$ {\tt -=} $\widehat{f}_{k+1}$} (GPsi1);

  \draw [connectline] (GPsi1) to [out=0,in=180] (PsiPsiForward);
  \draw [connectarrow] (PsiPsiForward) to [out=0,in=180] (GPsi2);

%
%
%
%

%
\end{scope}
\end{tikzpicture}

%% file: tikz/part1/psipsiConnect.tex
\begin{tikzpicture}
[fstyle/.style={shape=rectangle,rounded corners, 
                  text centered,
                  top color=\PsiModeColorTop,
                  bottom color=\PsiModeColorBottom,
                  draw=\PsiModeColorBorder},
gstyle/.style={shape=rectangle,rounded corners, 
                  text centered,
                  top color=\psiModeColorTop,
                  bottom color=\psiModeColorBottom,
                  draw=\psiModeColorBorder},
connectstyle/.style={shape=rectangle,rounded corners, 
                  text centered,
                  top color=\ConModeColorTop,
                  bottom color=\ConModeColorBottom,
                  text width=80pt,
                  draw=\ConModeColorBorder},
figstyle/.style={shape=rectangle,rounded corners, dashed,
                 text centered,
                 top color=\figrefColorTop,
                 bottom color=\figrefColorBottom,
                 minimum height = 80,
                 minimum width = 80,
                 text width=80pt,
                 draw=\figrefColorBorder},
sumstyle/.style={fill=blue},
connectarrow/.style={->,thick,>=triangle 45},
connectline/.style={-,thick},
connectdots/.style={-,very thick,dotted},
auto]

\def\PsiModeColorTop{red!40!white}
\def\PsiModeColorBottom{red!10!white}
\def\PsiModeColorBorder{red!40!black}
\def\ConModeColorTop{orange!40!white}
\def\ConModeColorBottom{orange!10!white}
\def\ConModeColorBorder{orange!40!black}
\def\psiModeColorTop{green!40!white}
\def\psiModeColorBottom{green!10!white}
\def\psiModeColorBorder{green!40!black}
\def\figrefColorTop{blue!10!white}
\def\figrefColorBottom{blue!10!white}
\def\figrefColorBorder{blue!90!white}

\def\braceOff{1.55}  

\def\fX{0}   
\def\gX{10} 
\def\fY{0}   
\def\gY{-6}   

\def\connectY{-3}  
\def\connectX{5}   

\def\PsiPsiFXOffset{4}  

\def\textY{-8.5}  
\def\textXL{1} 
\def\textXR{7} 

\begin{scope}[]

  \node at (\textXL,\textY) {$g \doteq f \times
  \left[\sqrt[\ast]{w^{(3)}}\right]^{-1}$};
  \node at (\textXR,\textY) {$h \doteq f \times
  \left[\sqrt[\ast]{w^{(5)}}\right]^{-1}$};
  \node (fpsi1) at (\fX,\fY) [fstyle] 
        {$\widehat{f}_k^{\,\psi,(3)}\equiv \widehat{g}_k^{\,\Psi,(3)}$};

  \node (PsiPsiBackward) at (\fX+\PsiPsiFXOffset,\fY) [figstyle] 
        {$\Psi$ - $\Psi$ connection Figure \ref{fig:fourier-connection}};

  \node (fpsi2) at (\fX+2*\PsiPsiFXOffset,\fY) [fstyle] 
        {$\widehat{g}_k^{\,\Psi,(0)}$};


  \node (connect) at (\connectX,\connectY) [connectstyle]
        {$\Psi^{(0)}_{}$ modes of $\left[\sqrt[\ast]{w^{(4)}}\right]^{-1} \times f$};


  \node (gpsi1) at (\gX-2*\PsiPsiFXOffset,\gY) [gstyle]
        {$\widehat{h}_k^{\,\Psi, (0)}$};

  \node (PsiPsiForward) at (\gX-\PsiPsiFXOffset,\gY) [figstyle] 
        {$\Psi$ - $\Psi$ connection\\ Figure \ref{fig:fourier-connection}};

  \node (gpsi2) at (\gX,\gY) [gstyle]
        {$\widehat{h}_k^{\,\Psi,(5)}\equiv \widehat{f}_k^{\,\psi,(5)}$};

  \draw [connectline] (fpsi1) to (PsiPsiBackward);
  \draw [connectarrow] (PsiPsiBackward) to (fpsi2);

  \draw [connectarrow] (fpsi2) to [out=-135,in=45] node 
        {$\widehat{g}_{k}$ {\tt -=} $\widehat{g}_{k+1}$} (connect);

  \draw [connectarrow] (connect) to [out=-135,in=45] node [swap]
        {$\widehat{g}_{k}$ {\tt -=} $\widehat{g}_{k+1}$} (gpsi1);

  \draw [connectline] (gpsi1) to (PsiPsiForward);
  \draw [connectarrow] (PsiPsiForward) to (gpsi2);
  
%
%

%
%

\end{scope}
\end{tikzpicture}

%% file: tikz/part1/GaussConstruct.tex
\begin{tikzpicture}[xscale=1.25]

\def\npi{-3.1415926535}
\newcommand\drawr[4]{\begin{scope}[xshift=#1,yshift=#2,xscale=#3,yscale=#4]
  \draw [line width=1pt] (-1,0) node [anchor=north] {{\scriptsize $\mathbf{-1}$}} --
        (1,0) node (lraxis) [anchor=north] {{\scriptsize $\mathbf{r=1}$}};
  \draw (-1,-3pt) -- (-1,3pt);
  \draw (1,-3pt) -- (1,3pt);
  \end{scope}
  }

\newcommand\drawt[4]{\begin{scope}[xshift=#1,yshift=#2,xscale=#3,yscale=#4]
  \draw [line width=1pt] (-1*\npi,0) node [anchor=north] {{\scriptsize $\mathbf{\pi}$}} --
        (1*\npi,0) node (lraxis) [anchor=north]
        {{\scriptsize $\mathbf{\theta=-\pi}$}};
  \draw (-\npi,-3pt) -- (-\npi,3pt);
  \draw (\npi,-3pt) -- (\npi,3pt);
  \end{scope}
  }

\def\gq{-0.44636196,  0.01488278,  0.45080338,  0.79036954, 0.97602054}

\def\grq{-0.50897301, -0.08716032,  0.32980906,  0.6822214 ,  0.9174362,  1.}

\def\teven{-2.03349201, -1.55591299, -1.10313117, -0.65938436, -0.21943523,
        0.21943523,  0.65938436,  1.10313117,  1.55591299,  2.03349201}

\def\todd{-2.1047876 , -1.65806739, -1.23469502, -0.81999974, -0.40920807,
        0, 0.40920807,  0.81999974,  1.23469502,  1.65806739,
        2.1047876}

\node (ne) at (2.5cm,-0.5cm) [anchor=north,
                              text centered,
                              rectangle,
                              rounded corners,
                              draw=blue!50!black,
                              top color=blue!30!white,
                              bottom color=blue!10!white]
                             {$N$ even};
\node at (ne.south) [anchor=north] {Jacobi-Gauss Quadrature};
\drawr{1cm}{-2cm}{1}{1}
\node (lrleft) at (1cm,-2.55cm) {};
\begin{scope}[xshift=1cm,yshift=-2cm]
  \foreach \t in \gq \draw [fill=blue] (\t,0) circle (2pt);
  \draw [decorate,decoration={brace,
                            mirror,
                            raise=0.6cm}] (-1cm,0cm) -- (1cm,0cm);
\end{scope}
\drawr{4cm}{-2cm}{-1}{1}
\node (lrright) at (4cm,-2.55cm) {};
\begin{scope}[xshift=4cm,yshift=-2cm,xscale=-1]
  \foreach \t in \gq \draw [fill=blue] (\t,0) circle (2pt);
  \draw [decorate,decoration={brace, 
                              raise=0.6cm}] (-1cm,0cm) -- (1cm,0cm);
\end{scope}

\node (rc) at (2.5cm,-4.2cm) {};
\draw [->] (lrleft) to [out=270,in=90] (rc);
\draw [->] (lrright) to [out=270,in=90] (rc);

\drawt{2.5cm}{-4.5cm}{0.8}{1}
\begin{scope}[xshift=2.5cm,yshift=-4.5cm,xscale=0.8]
  \foreach \t in \teven \draw [fill=blue] (\t,0) circle(2pt);
\end{scope}

\draw [dashed] (5.40cm,-0.5cm) -- (5.40cm,-6cm);
\begin{scope}[xshift=6cm]
  \node (ne) at (2.5cm,-0.5cm) [anchor=north,
                                text centered,
                                rectangle,
                                rounded corners,
                                draw=red!50!black,
                                top color=red!30!white,
                                bottom color=red!10!white]
                               {$N$ odd};
  \node at (ne.south) [anchor=north] {Jacobi-Gauss-Radau Quadrature};
  \drawr{1cm}{-2cm}{1}{1}
  \node (rrleft) at (1cm,-2.55cm) {};
  \begin{scope}[xshift=1cm,yshift=-2cm]
    \foreach \t in \grq \draw [fill=red] (\t,0) circle (2pt);
    \draw [decorate,decoration={brace,
                                mirror,
                                raise=0.6cm}] (-1cm,0cm) -- (1cm,0cm);
  \end{scope}
  \drawr{4cm}{-2cm}{-1}{1}
  \node (rrright) at (4cm,-2.55cm) {};
  \begin{scope}[xshift=4cm,yshift=-2cm,xscale=-1]
    \foreach \t in \grq \draw [fill=red] (\t,0) circle (2pt);
    \draw [decorate,decoration={brace, 
                                raise=0.6cm}] (-1cm,0cm) -- (1cm,0cm);
  \end{scope}
  \node (lc) at (2.5cm,-4.2cm) {};
  \draw [->] (rrleft) to [out=270,in=90] (lc);
  \draw [->] (rrright) to [out=270,in=90] (lc);

  \drawt{2.5cm}{-4.5cm}{0.8}{1}
  \begin{scope}[xshift=2.5cm,yshift=-4.5cm,xscale=0.8]
    \foreach \t in \todd \draw [fill=red] (\t,0) circle(2pt);
  \end{scope}
\end{scope}

\end{tikzpicture}

%% file: tikz/part1/quadpoints.tex
\newcommand\pgfmathsinandcos[3]{%
  \pgfmathsetmacro#1{sin(#3)}%
  \pgfmathsetmacro#2{cos(#3)}%
}

\newcommand\drawcircle[1]{
  \draw [dashed] (0,0) node [anchor=center] {#1} circle (1cm);
  }

\newcommand\drawaxis[1]{
  \draw [<->,color=gray] (-15*\xscale,0) -- (15*\xscale,0);
  \node at (-15*\xscale,0) [anchor=east] {#1};
  }

\begin{tikzpicture}
\def\sina{0}
\def\cosa{0}
\def\radang{57.295779513082323}
\def\xscale{0.45}

\def\thetazero{
       -2.991993 , -2.6927937, -2.3935944, -2.0943951, -1.7951958,
       -1.4959965, -1.1967972, -0.8975979, -0.5983986, -0.2991993,
        0.       ,  0.2991993,  0.5983986,  0.8975979,  1.1967972,
        1.4959965,  1.7951958,  2.0943951,  2.3935944,  2.6927937,
        2.991993 }
\def\xzero{
        -13.34407264,  -4.38128627,  -2.54795845,  -1.73205081,
        -1.25396034,  -0.9278644 ,  -0.68178846,  -0.48157462,
        -0.30845915,  -0.15072575,   0.        ,   0.15072575,
         0.30845915,   0.48157462,   0.68178846,   0.9278644 ,
         1.25396034,   1.73205081,   2.54795845,   4.38128627,
         13.34407264}

\def\thetatwo{
        -2.75061425, -2.46940457, -2.19280076, -1.91765367, -1.64315354,
       -1.3689941 , -1.09503242, -0.82119161, -0.54742536, -0.27370261,
        0.        ,  0.27370261,  0.54742536,  0.82119161,  1.09503242,
        1.3689941 ,  1.64315354,  1.91765367,  2.19280076,  2.46940457,
        2.75061425}
\def\xtwo{
        -5.05004233,-2.86247354, -1.94738729, -1.42479679, -1.07510726,
        -0.81612639,
         -0.60969296, -0.43533974, -0.28075932, -0.13771209,  0.0,
         0.13771209,
          0.28075932,  0.43533974,  0.60969296,  0.81612639,
          1.07510726,  1.42479679,
          1.94738729,  2.86247354, 5.05004233}

\def\thetafour{
        -2.58074662, -2.30550504, -2.04214504, -1.78319802, -1.5264179 ,
       -1.27085944, -1.0160452 , -0.76170184, -0.50765607, -0.25378682,
        0.        ,  0.25378682,  0.50765607,  0.76170184,  1.0160452 ,
        1.27085944,  1.5264179 ,  1.78319802,  2.04214504,  2.30550504,
        2.58074662}
\def\xfour{
        -3.47207303, -2.25109488, -1.63205282, -1.23864373, -0.95657795,
       -0.73746407, -0.55676549, -0.40039972, -0.25942355, -0.1275789 ,
        0.        ,  0.1275789 ,  0.25942355,  0.40039972,  0.55676549,
        0.73746407,  0.95657795,  1.23864373,  1.63205282,  2.25109488,
        3.47207303}

\def\thetasix{
        -2.44370032, -2.17430824, -1.92147158, -1.67527223, -1.43251158,
       -1.19176755, -0.95228868, -0.713628  , -0.47549268, -0.23767233,
        0.        ,  0.23767233,  0.47549268,  0.713628  ,  0.95228868,
        1.19176755,  1.43251158,  1.67527223,  1.92147158,  2.17430824,
        2.44370032}
\def\xsix{
        -2.7485009 , -1.9038588 , -1.4305968 , -1.11034023, -0.87046511,
       -0.67811095, -0.51571995, -0.3727698 , -0.24232939, -0.11939875,
        0.        ,  0.11939875,  0.24232939,  0.3727698 ,  0.51571995,
        0.67811095,  0.87046511,  1.11034023,  1.4305968 ,  1.9038588 ,
        2.7485009}

\def\thetaeight{
        -2.32839851, -2.06499962, -1.82123585, -1.58570651, -1.35458817,
       -1.1261236 , -0.89935693, -0.67370486, -0.4487764 , -0.22428504,
        0.        ,  0.22428504,  0.4487764 ,  0.67370486,  0.89935693,
        1.1261236 ,  1.35458817,  1.58570651,  1.82123585,  2.06499962,
        2.3283985}
\def\xeight{
        -2.32238725, -1.67471531, -1.28801112, -1.01502245, -0.80419574,
       -0.63122304, -0.48265857, -0.35019938, -0.22823161, -0.11261499,
        0.        ,  0.11261499,  0.22823161,  0.35019938,  0.48265857,
        0.63122304,  0.80419574,  1.01502245,  1.28801112,  1.67471531,
        2.32238725}

\def\czero{black}
\def\ctwo{blue}
\def\cfour{green}
\def\csix{red}
\def\ceight{yellow}

\begin{scope}[xshift=-6cm]
  \drawcircle{$\gamma=0$}
  \foreach \t in \thetazero {
      \pgfmathsinandcos\sina\cosa{\t*\radang+90}
      \draw [fill=\czero] (\sina,\cosa) circle (2pt);
  }

  \begin{scope}[xshift=3cm]
    \drawcircle{$\gamma=2$}
    \foreach \t in \thetatwo {
        \pgfmathsinandcos\sina\cosa{\t*\radang+90}
        \draw [fill=\ctwo] (\sina,\cosa) circle (2pt);
    }
  \end{scope}

  \begin{scope}[xshift=6cm]
    \drawcircle{$\gamma=4$}
    \foreach \t in \thetafour {
        \pgfmathsinandcos\sina\cosa{\t*\radang+90}
        \draw [fill=\cfour] (\sina,\cosa) circle (2pt);
    }
  \end{scope}
  \begin{scope}[xshift=9cm]
    \drawcircle{$\gamma=6$}
    \foreach \t in \thetasix {
        \pgfmathsinandcos\sina\cosa{\t*\radang+90}
        \draw [fill=\csix] (\sina,\cosa) circle (2pt);
    }
  \end{scope}
  \begin{scope}[xshift=12cm]
    \drawcircle{$\gamma=8$}
    \foreach \t in \thetaeight {
        \pgfmathsinandcos\sina\cosa{\t*\radang+90}
        \draw [fill=\ceight] (\sina,\cosa) circle (2pt);
    }
  \end{scope}
\end{scope}

\begin{scope}[yshift=-2cm]
  \drawaxis{$s=1$}
  \foreach \x in \xzero {
    \draw [fill=\czero] (\x*\xscale,0) circle(2pt);
  }

  \begin{scope}[yshift=-1cm] 
    \drawaxis{$s=3$}
    \foreach \x in \xtwo {
      \draw [fill=\ctwo] (\x*\xscale,0) circle(2pt);
    }
  \end{scope}

  \begin{scope}[yshift=-2cm] 
    \drawaxis{$s=5$}
    \foreach \x in \xfour {
      \draw [fill=\cfour] (\x*\xscale,0) circle(2pt);
    }
  \end{scope}

  \begin{scope}[yshift=-3cm] 
    \drawaxis{$s=7$}
    \foreach \x in \xsix {
      \draw [fill=\csix] (\x*\xscale,0) circle(2pt);
    }
  \end{scope}
  
  \begin{scope}[yshift=-4cm] 
    \drawaxis{$s=9$}
    \foreach \x in \xeight {
      \draw [fill=\ceight] (\x*\xscale,0) circle(2pt);
    }
  \end{scope}
  \begin{scope}[yshift=3.5cm]
    \node[anchor=south] at (0,0) {$s \doteq \gamma + 1$};
  \end{scope}
\end{scope}

\end{tikzpicture}

%% file: part1Content/semiinfinite.tex
\label{sec:semiinfinite}
The generalized Wiener basis functions we have derived can be used for function
expansions on the infinite line. In order to address expansions on semi-infinite
intervals, we can instead use either the even or odd Jacobi polynomial basis
sets that make up the Fourier functions constructed in Section
\ref{sec:derivationFourier}. 

The Jacobi functions from Lemma \ref{lemma:jacobi-functions} can be mapped and
weighted in a procedure identical to the construction of the Wiener basis. The
result is the collection of functions 

\begin{equation}
\label{eq:rhodef}
\begin{array}{rclr}
  \rho_n^{(s)} & = &\sqrt{w_x^{(s)}} \Psi_n^{(s)}(x) & \\[8pt]
              & = &\left(\frac{2}{x^2+1}\right)^{s/2} \tilde{P}_n^{(-1/2,s-3/2)}
                    \left(\frac{1-x^2}{1+x^2}\right), & n \in \mathbb{N}_0
\end{array}
\end{equation}

\noindent These functions are a direct mapping and weighting of the Jacobi polynomials.
Because of this, they are orthonormal and complete in
$L^2\left(\mathbb{R}^+,\mathbb{R}\right)$. Mapping techniques for classical
functions are not novel and we discuss existing methods in Section
\ref{sec:mjpoly}. The classical competitor for spectral expansions on
semi-infinite intervals is the set of Laguerre functions (weighted Laguerre
polynomials). A comparison between the Laguerre functions and the functions
defined in (\ref{eq:rhodef}) will be made in Part \Rmnum{2}, and in Section
\ref{sec:mjpoly} a different mapping transformaing Jacobi polynomials to the
semi-infinite line will be addressed. 

We make use of the regular square root function $\sqrt{w_x^{(s)}}$ in
(\ref{eq:rhodef}) instead of the phase-shifted version
$\sqrt[\asterisk]{w_x^{(s)}}$ because there is no need to have complex-valued
functions. The phase-shifted square root was a convenient choice for the Wiener
functions on the infinite line: its compact Fourier series representation
(\ref{eq:phase-shift-rewrite}) enabled fast connections (Section
\ref{sec:propjacobi-connections}) and sparse differentation matrices (Section
\ref{sec:propjacobi-stiffness}). By using the real-valued square root in
(\ref{eq:rhodef}) we sacrifice these two properties. However, the FFT can still
be used for the evaluation of modal coefficients if $s$ is an integer. 

The caveat in using these functions for expansions on the semi-infinite interval
is the fact that they all have zero-valued odd derivatives at $x=0$. This
parallels the same property at $\theta=0$ for a cosine series on $\theta\in
[0,\pi]$. Alternative mappings of the Jacobi polynomials to the semi-infinite
line do not exhibit this restriction, but those mappings also preserve the
$\mathcal{O}(N^2)$ time-stepping restriction for nodal-based polynomial solvers
of time-dependent partial differential equations using explicit time-integration
on finite intervals. In constrast, the functions (\ref{eq:rhodef}) only have an
$\mathcal{O}(N)$ time-step restriction, similar to the time-step restriction for
a finite-interval cosine basis expansion. 

The restriction of the Wiener functions to the semi-infinite interval as defined
in (\ref{eq:rhodef}) comes both with advantages and sacrifices. These functions
are purely weighted maps of Jacobi polynomials and are therefore easy to
implement. Some of the attractive features of the Wiener rational basis
functions on infinite intervals are lost (e.g. sparse stiffness matrices).
However, these functions have properties that are advantageous when compared
with existing mapping techniques (Section \ref{sec:mjpoly}). A numerical
comparison between those mapping techniques, the functions (\ref{eq:rhodef}),
and the Laguerre functions will be made in Part \Rmnum{2}.

%% file: part1Content/mjpoly.tex
\label{sec:mjpoly}
Before concluding this article with a summary of the derived properties of the
generalized Wiener basis, we first summarize existing results on the topic of
mapping Jacobi polynomials from the finite interval to the infinite interval.
This method is very closely related to our strategy of mapping a generalized
Fourier series from the canonical finite Fourier interval to the real line.
Numerical studies comparing these methods are presented in Part \Rmnum{2}, but it is
appropriate to acknowledge these functions here, and to discuss how they relate to
the Wiener rational function basis. 

\subsection{The Infinite Interval}
\label{sec:mjpoly-infinite}
The main idea for our generalization of Wiener's original rational basis is
using a `well-behaved' mapping to transform functions on a finite interval to
those on an infinite interval. This basic idea is classical
{\cite{gottlieb1977}}. Indeed one of the more popular mappings that has gained
momentum in the literature are the so-called `mapped Chebyshev'
functions/polynomials.

In order to further generalize the mapped Chebyshev functions, we will briefly
restate their derivation from our point of view. We begin with the Jacobi
polynomials $P_n^{(\alpha, \beta)} (r)$ on $r \in [- 1, 1]$. Mapping via $r =
\cos \theta$ to $\theta \in [0, \pi]$ yields trigonometric polynomials. We now
`stretch' the domain to $\Theta \in [- \pi, \pi]$ via the affine mapping $\Theta = 2
\theta - \pi$. Finally, we utilize the usual linear fractional map $e^{i
\Theta} = \frac{i - x}{i + x}$ (i.e. rotation of the Riemann Sphere) to yield
functions on the real line $x \in \mathbbm{R}$. For all $s, t > \frac{1}{2}$,
this results in the functions $\tmop{PB}_n^{(s, t)} (x)$, defined as
\[ \tmop{PB}_n^{(s, t)} (x) = \widetilde{P}_n^{((2 s - 3) / 2, (2 t - 3) / 2)} \left(
   \frac{x}{\sqrt{1 + x^2}} \right), \]
orthonormal on the real line under the weight
\[ w_{\tmop{PB}}^{(s, t)} = \left[ 1 - \frac{x}{\sqrt{1 + x^2}} \right]^{(2 s
   - 3) / 2} \left[ 1 + \frac{x}{\sqrt{1 + x^2}} \right]^{(2 t - 3) / 2}, \]
and the weighted functions
\[ \tmop{pb}_n^{(s, t)} \assign \sqrt{w_{\tmop{PB}}^{(s, t)}}
   \tmop{PB}_n^{(s, t)}, \]
are orthonormal under the unweighted inner product. When $s = t = 1$, the
functions $\tmop{PB}_n^{(s, t)}$ coincide with the mapped Chebyshev polynomials
$\tmop{TB}_n (x)$ introduced in {\cite{boyd1982}} and subsequently developed
in {\cite{boyd1987c}} and {\cite{boyd1989}}, although the original idea of
applying spectral expansions over finite intervals to solving problems over
infinite intervals seems to come from {\cite{grosch1977}}. In any case, the
mapped Jacobi functions $\tmop{pb}^{(s, t)}_n$ decay like $\frac{1}{|x|^s}$
for $x \rightarrow - \infty$ and $\frac{1}{|x|^t}$ for $x \rightarrow +
\infty$. The advantage of these functions is that the decay can be different
as $|x| \rightarrow \infty$. Also, others have already explored some convergence theory in
function spaces {\cite{benyu2002}} and applications to differential equations
{\cite{zhongqing2002}} for the Chebyshev case $s = t = 1$. In Part \Rmnum{2} when we
present numerical examples, we will use the basis set $\tmop{pb}_n^{(s, t)}$
with $s = t = 1$, i.e. the Chebyshev case.

Note that because all of these mapped types of polynomials and the generalized
Wiener basis we have presented ultimately stem from Jacobi polynomials and
mappings of similar character, all these basis sets are related in some
fashion. To relate the mapped Jacobi functions to the generalized Wiener
rational functions, we have
\[ \tmop{PB}_n^{(s, s)} (x) \propto \tmop{Re} \left\{ \Phi_n^{(s)} \left(
   \frac{x + \sqrt{x^2 + 1} - 1}{x - \sqrt{x^2 + 1} + 1} \right) \right\} \]
In Table \ref{tab:function-relations} we relate the unweighted functions to
the generalized Wiener rational basis, modulo multiplicative constants. In
this article we make no observations about how mapped Jacobi polynomials
compare to the Wiener basis set as a practical tool for function expansions.
However, such a comparison will be a central theme in Part \Rmnum{2}.

\begin{table}[tbp]
 \renewcommand{\arraystretch}{1.8}
  \begin{tabular}{|c|c|c|c|c|}\hline
    Previous function & Name/classification & Interval & Reference &
    Relation\\\hline
    $T \nonesep B_n$ & Cheyshev rational functions (1st)  & $\mathbbm{R}$ &
    {\cite{boyd1990a}}, {\cite{boyd1982}}, {\cite{boyd1987}} &
    $\tmop{PB}^{(1, 1)}_n$\\\hline
    $S \nonesep B_n$/$U \nonesep B_n$ & Chebyshev rational functions (2nd)
    & $\mathbbm{R}$ & {\cite{boyd1990a}}, {\cite{boyd1987}}, {\cite{boyd1989}}
    & $\tmop{PB}^{(2, 2)}_n$\\\hline
    $C_n$/$\tmop{CC}_n$ & Christov functions (even) & $\mathbbm{R}$, $[0,
    \infty)$ & {\cite{christov1982}}, {\cite{boyd1990a}} & $\tmop{Im} \left\{
    \phi^{(1, 0)}_n \right\}$\\\hline
    $S_n / \tmop{SC}_n$ & Christov functions (odd) & $\mathbbm{R}$, $[0,
    \infty)$ & {\cite{christov1982}}, {\cite{boyd1990a}} & $\tmop{Re} \left\{
    \phi^{(1, 0)}_n \right\}$\\\hline
    $C \nonesep H_n$ & Higgins functions (even) & $\mathbbm{R}$, $[0,
    \infty)$ & {\cite{boyd1990a}} & $\tmop{Re} \left\{ \Phi^{(1, 0)}_n
    \right\}$\\\hline
    $S \nonesep H_n$ & Higgins functions (odd) & $\mathbbm{R}$, $[0,
    \infty)$ & {\cite{boyd1990a}} & $\tmop{Im} \left\{ \Phi^{(1, 0)}_n
    \right\}$\\\hline
    $\rho_k$ & (Complex) Higgins functions & $\mathbbm{R}$ &
    {\cite{higgins1977}}, {\cite{christov1982}} & $\Phi^{(1,
    0)}_k$\\\hline
    $\sigma_k$ & (Complex) Wiener rational functions & $\mathbbm{R}$ &
    {\cite{wiener1949}}, {\cite{christov1982}} & $\phi^{(1,
    0)}_k$\\\hline
    $T \nonesep L_n$ & Half-infinte Chebyshev rational functions & $[0,
    \infty)$ & {\cite{boyd1987b}} & $\text{PL}_n^{(1/2)}$ \\\hline
  \end{tabular}
  \caption{Relationship between orthogonal functions in previous work and the
  current bases presented.\label{tab:function-relations}}
\end{table}

\subsection{The Semi-Infinite Interval}
\label{sec:mjpoly-semiinfinite}
To perform spectral expansions on semi-infinite intervals, the only classical
technique is the Laguerre polynomial/function method. However, mapping
techniques can be used to transform finite-interval methods to semi-infinite
interval methods. 

As with Section \ref{sec:mjpoly-infinite}, we explain the choice of mapping from
our point of view as a mapping of the Riemann Sphere. The Jacobi polynomials are
defined on $r\in [-1,1]$. If we allow complex values of $r$, then we may consider using a
linear fractional map to transform the Jacobi polynomial domain to the
semi-infinite line. The ordered assignments $r = \{1,0,-1\}$ to $x=\{0, 1,
\infty\}$ specify the transformation uniquely as 

\begin{align}
\label{eq:semiinfinite-boyd-map}
\begin{array}{ccc}
x = \frac{1-r}{1+r} &  & r = \frac{1-x}{1+x}.
\end{array}
\end{align}

\noindent If necessary, one can also specify the relationship to $\theta$ and the cosine
series on $[0,\pi]$. For details, see \cite{boyd1987b}. Our definition of the
transformation differs only in orientiation from that presented in
\cite{boyd1987b}. We have chosen this orientation so that the Jacobi parameter
$\beta$ is assigned to the location $x=\infty$ in order to mimic to the same
assignment for the Wiener functions. 

In the literature the maps of the Chebyshev polynomials under the transformation
(\ref{eq:semiinfinite-boyd-map}) are labeled $\text{TL}_n(x)$. Adopting similar
notation, we define 
\begin{align*}
\text{PL}_n^{(s)}(x) = \tilde{P}_n^{(-1/2,2 s-2)}\left(\frac{1-x}{1+x}\right),
\hspace{1cm} x\in[0,\infty],
\end{align*}

\noindent which are $L^2$-complete and orthonormal under the weight function 
\begin{align*}
  w_{\text{PL}}^{(s)}(x) = \frac{1}{2\sqrt{x}}\left(\frac{2}{1+x}\right)^{(2s)}. 
\end{align*}

\noindent It is then possible to define the weighted functions 
\begin{align}
\label{eq:pldef}
\text{pl}_n^{(s)}(x) &= \left(\frac{2}{1+x}\right)^s \text{PL}_n^{(s)}(x) \\
  &= \left(\frac{2}{1+x}\right)^s
     \tilde{P}_n^{(-1/2,2 s-2)}\left(\frac{1-x}{1+x}\right),
\end{align}

\noindent which are $L^2$-complete and orthonormal under the weighted $L^2$ inner
product 
\begin{align*}
\left\langle f,g \right\rangle_{w_{PL}^{(0)}} = \int_0^\infty f\,g\,\frac{1}{2\sqrt{x}}
  \mathrm{d}x,
\end{align*}

\noindent for any $s>\frac{1}{2}$. The $\text{pl}_n^{(s)}$ are defined for $x\in [0,\infty]$ and decay
like $x^s$ as $x\rightarrow \infty$. A significant difference between the
Wiener-type functions (both on the infinite and semi-infinte intervals) and the
set defined in (\ref{eq:pldef}) is the fact that these functions are not
orthogonal in the unweighted $L^2$ inner product, but instead in the 
norm defined by the above inner product. This choice was made (as opposed to
defining functions in the unweighted inner product) to ensure that
integer values of $s$ resulted in a Jacobi polynomial family that was amenable
to usage of the FFT. 

The main observation we make regarding this basis is that these functions are
the result of a linear fractional map directly from the Jacobi domain;
therefore, they will inherit the $\mathcal{O}(N^2)$ time-step restriction of
nodal explicit time-integration methods for time-dependent partial differential
equations. The same observation can be made about the functions defined in
\cite{boyd1987b}.

%% file: part1Content/conclusion.tex
\label{sec:conclusion}We have presented a collection of generalized Fourier
series which, when mapped and weighted appropriately, generates a basis set on
the infinite interval with a tunable rate of decay. For each rate of decay $s$
satisfying $s > \frac{1}{2}$ the resulting basis set $\phi_k^{(s)}$:
\begin{itemize}
  \item is orthonormal and complete in $L^2 \left( \mathbbm{R}, \mathbbm{C}
  \right)$
  
  \item is characterized by $x^{- s}$ decay for $|x| \rightarrow \infty$
  
  \item can be generated via Jacobi polynomial recurrence relations
  
  \item has sparse connection properties that can be efficiently exploited via
  combinations of sparse Fourier and Jacobi connections
  
  \item has an $N \times N$ Galerkin stiffness/differentiation matrix that has
  at most $6 N$ nonzero entries with $\mathcal{O} (N)$ spectral radius
  
  \item is characterized by a `Gauss-like' quadrature rule.
\end{itemize}
When $s \in \mathbbm{N}$, the basis set is a rational function; we will show
in Part \Rmnum{2} that in this case we can use the FFT for modal-nodal
transformations. The case $s = 1$ corresponds to a mapping and weighting of
the canonical Fourier series, as discovered by others previously. Due to the
original presentation of the $s = 1$ basis by Wiener {\cite{wiener1949}}, we
call the functions $\phi_k^{(s)}$ the generalized Wiener rational basis
functions.

These basis functions have a similar flavor to directly mapped and weighted
Jacobi polynomials (called $\tmop{pb}_n^{(s, t)}$ here). In Part \Rmnum{2} we will
compare these basis sets and discuss advantages and disadvantages of each. In
addition, we will also employ the Sinc and Hermite functions in test cases in an
attempt to investigate a relatively broad class of spectral approximation
methods. In contrast to \cite{shen2009} which reviews much of the theory
present for expansions on the infinite interval, we concentrate on numerical issues,
including application of the FFT. We will extend our investigation to the
semi-infinite interval to compare the Laguerre polynomials/functions, the mapped
Jacobi functions (denoted $\text{pl}_n^{(s)}$ here), and the restriction of the
Wiener functions to the semi-infinite interval as given in Section
\ref{sec:semiinfinite}.

We do not wish to claim that, on the infinite or semi-infinite intervals,
genuinely global spectral expansions are truly superior to alternative numerical
approximations; rather we wish to identify the generalized Wiener basis set as a
novel competitor to existing global spectral expansions. Part \Rmnum{2} will
follow up to show that the Wiener basis set is very competitive with existing
expansions.

\subsubsection*{Acknowledgements}
The authors acknowledge partial support for this work by
AFOSR award FA9550-07-1-0422.

%% file: part1Content/appendix.tex
\section{Recurrence Coefficients}

\label{app:recurrence}In this appendix we compile various recurrence relations
for the Jacobi/Szeg$\ddot{\text{o}}$-Fourier/Wiener rational functions. We state the
recurrences in terms of the Szeg$\ddot{\text{o}}$-Fourier functions
$\Psi_k^{(\gamma)}$, but note that they all apply equally well to the
unweighted Wiener rational functions as well. Note that we only list
recurrences for $k \geq 0$; for $k < 0$, we may use the conjugation relation
(\ref{eq:indsym-1}) to obtain $\Psi_{- |k|}^{(\gamma)}$ at almost no
additional cost. We first require a tour of some Jacobi polynomials
recurrences:

\begin{eqnarray}
  \sqrt{b_{n + 1}^{(\alpha, \beta)}} \tilde{P}_{n + 1}^{(\alpha, \beta)} & = &
  \left[ r - a_n^{(\alpha, \beta)} \right] \tilde{P}_n^{(\alpha, \beta)} -
  \sqrt{b_n^{(\alpha, \beta)}} \tilde{P}_{n - 1}^{(\alpha, \beta)},
  \label{eq:opoly-3term}\\
  &  &  \nonumber\\
  (1 - r^2) \tilde{P}_n^{(\alpha, \beta)} & = & \sum_{i = 0}^2 \varepsilon_{n,
  i}^{(\alpha, \beta)} \tilde{P}_{n + i}^{(\alpha - 1, \beta - 1)},
  \label{eq:opoly-demote}\\
  &  &  \nonumber\\
  \tilde{P}_n^{(\alpha, \beta)} & = & \sum_{i = 0}^2 \eta_{n, - i}^{(\alpha,
  \beta)} \tilde{P}_{n - i}^{(\alpha + 1, \beta + 1)},
  \label{eq:opoly-promote} \\
  & & \nonumber\\
  (1 - r) \tilde{P}_n^{(\alpha, \beta)} & = & \mu_{n, 0}^{(\alpha, \beta)} 
  \tilde{P}_n^{(\alpha - 1, \beta)} - \mu_{n, 1}^{(\alpha, \beta)}
  \tilde{P}_{n + 1}^{(\alpha - 1, \beta)},  \label{eq:jdemotiona}\\
  &  &  \nonumber\\
  (1 + r) \tilde{P}_n^{(\alpha, \beta)} & = & \mu_{n, 0}^{(\beta, \alpha)} 
  \tilde{P}_n^{(\alpha, \beta - 1)} + \mu_{n, 1}^{(\beta, \alpha)}
  \tilde{P}_{n + 1}^{(\alpha, \beta - 1)},  \label{eq:jdemotionb}\\
  &  &  \nonumber\\
  \tilde{P}_n^{(\alpha, \beta)} & = & \nu_{n, 0}^{(\alpha, \beta)}
  \tilde{P}^{(\alpha + 1, \beta)}_n - \nu_{n, - 1}^{(\alpha, \beta)}
  \tilde{P}_{n - 1}^{(\alpha + 1, \beta)},  \label{eq:jpromotiona}\\
  &  &  \nonumber\\
  \tilde{P}_n^{(\alpha, \beta)} & = & \nu_{n, 0}^{(\beta, \alpha)}
  \tilde{P}_n^{(\alpha, \beta + 1)} + \nu_{n, - 1}^{(\beta, \alpha)}
  \tilde{P}_{n - 1}^{(\alpha, \beta + 1)},  \label{eq:jpromotionb}\\
  &  &  \nonumber\\
  \frac{\mathd}{\mathd r} \tilde{P}^{(\alpha, \beta)}_n & = &
  \gamma_n^{(\alpha, \beta)} \tilde{P}_{n - 1}^{(\alpha + 1, \beta + 1)}, 
  \label{eq:jdiff}
\end{eqnarray}
where $\mu_{n, 0 / 1}^{(\alpha, \beta)}$, $\nu_{n, 0 / - 1}^{(\alpha,
\beta)}$, and $\gamma_n^{(\alpha, \beta)}$ in 
(\ref{eq:jdemotiona})-(\ref{eq:jdiff}) are constants for which we take explicit
formulae from {\cite{narayan2009}}:

\begin{eqnarray}
  \mu_{n, 0}^{(\alpha, \beta)} & = & \sqrt{\frac{2 (n + \alpha) (n + \alpha +
  \beta)}{(2 n + \alpha + \beta) (2 n + \alpha + \beta + 1)}},
  \label{eq:mu0}\\
  &  &  \nonumber\\
  \mu_{n, 1}^{(\alpha, \beta)} & = & \sqrt{\frac{2 (n + 1) (n + \beta + 1)}{(2
  n + \alpha + \beta + 1) (2 n + \alpha + \beta + 2)}},  \label{eq:mu1}\\
  &  &  \nonumber\\
  \nu_{n, 0}^{(\alpha, \beta)} & = & \sqrt{\frac{2 (n + \alpha + 1) (n +
  \alpha + \beta + 1)}{(2 n + \alpha + \beta + 1) (2 n + \alpha + \beta + 2)}},
  \label{eq:nu0}\\
  &  &  \nonumber\\
  \nu_{n, - 1}^{(\alpha, \beta)} & = & \sqrt{\frac{2 n (n + \beta)}{(2 n +
  \alpha + \beta) (2 n + \alpha + \beta + 1)}},  \label{eq:nu1}\\
  &  &  \nonumber\\
  \gamma_n^{(\alpha, \beta)} & = & \sqrt{n (n + \alpha + \beta + 1)}. 
  \label{eq:eta}
\end{eqnarray}
The three-term recurrence coefficients in (\ref{eq:opoly-3term}) are given by
{\cite{gautschi2004}}:
\begin{eqnarray}
  a_n^{(\alpha, \beta)} & = & \left\{ \begin{array}{lll}
    \frac{\beta - \alpha}{\alpha + \beta + 2}, &  & n = 0,\\
    &  & \\
    \frac{\beta^2 - \alpha^2}{(2 n + \alpha + \beta) (2 n + \alpha + \beta +
    2)}, &  & n > 0.
  \end{array}  \label{eq:opoly-an} \right.\\
  &  &  \nonumber\\
  b_n^{(\alpha, \beta)} & = & \left\{ \begin{array}{lll}
    \frac{2^{\alpha + \beta + 1} \Gamma (\alpha + 1) \Gamma (\beta +
    1)}{\Gamma (\alpha + \beta + 2)}, &  & n = 0,\\
    &  & \\
    \frac{4 (\alpha + 1) (\beta + 1)}{(\alpha + \beta + 2)^2 (\alpha + \beta +
    3)}, &  & n = 1,\\
    &  & \\
    \frac{4 n (n + \alpha) (n + \beta) (n + \alpha + \beta)}{(2 n + \alpha +
    \beta - 1) (2 n + \alpha + \beta)^2 (2 n + \alpha + \beta + 1)}, &  & n >
    1.
  \end{array}  \label{eq:opoly-bn} \right.
\end{eqnarray}
The demotion recurrence coefficients in (\ref{eq:opoly-demote}) can be
obtained by determining the analogous \ relations for the monic orthogonal
polynomials ({\cite{abramowitz1972}}, {\cite{szego1975}}) and then employing
the appropriate normalizations:

\begin{eqnarray}
  \varepsilon_{n, 0}^{(\alpha, \beta)} & = & \left\{ \begin{array}{lll}
    2 \sqrt{\frac{\alpha \beta}{(\alpha + \beta) (\alpha + \beta + 1)}}, &  &
    n = 0,\\
    &  & \\
    \frac{2}{(\alpha + \beta + 2)} \sqrt{\frac{(\alpha + 1) (\beta + 1)
    (\alpha + \beta)}{(\alpha + \beta + 3)}}, &  & n = 1,\\
    &  & \\
    \frac{2}{(2 n + \alpha + \beta)} \sqrt{\frac{(n + \alpha) (n + \beta) (n +
    \alpha + \beta - 1) (n + \alpha + \beta)}{(2 n + \alpha + \beta - 1) (2 n
    + \alpha + \beta + 1)}}, &  & n > 1.
  \end{array}  \right.\\
  &  &  \nonumber\\
  \varepsilon_{n, 1}^{(\alpha, \beta)} & = & \left\{ \begin{array}{lll}
    \frac{2 (\alpha - \beta)}{(\alpha + \beta + 2) \sqrt{\alpha + \beta}}, &  &
    n = 0,\\
    &  & \\
    \frac{2 (\alpha - \beta) \sqrt{(n + 1) (n + \alpha + \beta)}}{(2 n +
    \alpha + \beta) (2 n + \alpha + \beta + 2)}, &  & n > 0.
  \end{array}  \right.\\
  &  &  \nonumber\\
  \varepsilon_{n, 2}^{(\alpha, \beta)} & = & \left\{ \begin{array}{lll}
    \frac{2}{\alpha + \beta + 2} \sqrt{\frac{2 (\alpha + 1) (\beta +
    1)}{(\alpha + \beta + 1) (\alpha + \beta + 3)}}, &  & n = 0,\\
    &  & \\
    \frac{2}{2 n + \alpha + \beta + 2} \sqrt{\frac{(n + 1) (n + 2) (n + \alpha
    + 1) (n + \beta + 1)}{(2 n + \alpha + \beta + 1) (2 n + \alpha + \beta +
    3)}}, &  & n > 0.
  \end{array}  \right.
\end{eqnarray}
Finally, the promotion relation (\ref{eq:opoly-promote}) coefficients can also
be determined:
\begin{eqnarray*}
  \eta_{n, 0}^{(\alpha, \beta)} & = & \varepsilon_{n, 0}^{(\alpha + 1, \beta +
  1)},\\
  &  & \\
  \eta_{n, - 1}^{(\alpha, \beta)} & = & \varepsilon_{n - 1, 1}^{(\alpha + 1,
  \beta + 1)},\\
  &  & \\
  \eta_{n, - 2}^{(\alpha, \beta)} & = & - \varepsilon_{n - 2, 2}^{(\alpha + 1,
  \beta + 1)} .
\end{eqnarray*}
Of course, (\ref{eq:opoly-demote})-(\ref{eq:opoly-promote}) are
consequences of combining 
(\ref{eq:jdemotiona})-(\ref{eq:jpromotionb}). Using the orthogonal polynomial
three-term recurrence relation (\ref{eq:opoly-3term}) we can show the
following recurrence relation for the Szeg$\ddot{\text{o}}$-Fourier functions
$\Psi_n^{(\gamma, \delta)} (\theta)$:

\begin{equation}
  \label{eq:Psi-cosrecurrence} \begin{array}{lll}
    \Psi_{n + 1}^{(\gamma)} & = & \left[ U_n^{(\gamma)} \cos \theta -
    V_n^{(\gamma)} \right] \Psi_n^{(\gamma)} + \left[ U_{- n}^{(\gamma)} \cos
    \theta - V_{- n}^{(\gamma)} \right] \Psi_{- n}^{(\gamma)}\\
    &  & \\
    &  & - W_n^{(\gamma)} \Psi_{n - 1}^{(\gamma)} - W_{- n}^{(\gamma)}
    \Psi_{- (n - 1)}^{(\gamma)} .
  \end{array}
\end{equation}
In the following expressions, we make use of the following definitions: for a
given $\gamma > - \frac{1}{2}$,
\[ \begin{array}{lll}
     \alpha \assign - \frac{1}{2}, &  & \beta \assign \gamma - \frac{1}{2} .
   \end{array} \]
The recurrence constants are then given by
\[ \begin{array}{lll}
     U_{\pm n}^{(\gamma)} & = & \frac{1}{2} \left[ \sqrt{\frac{1}{b_{n +
     1}^{(\alpha, \beta)}}} \pm \sqrt{\frac{1}{b_n^{(\alpha + 1, \beta + 1)}}}
     \right],\\
     &  & \\
     V_{\pm n}^{(\gamma)} & = & \pm \frac{1}{2} \left[ \frac{a_n^{(\alpha,
     \beta)}}{\sqrt{b_{n + 1}^{(\alpha, \beta)}}} + \frac{a_{n - 1}^{(\alpha +
     1, \beta + 1)}}{\sqrt{b_n^{(\alpha + 1, \beta + 1)}}} \right],\\
     &  & \\
     W_{\pm n}^{(\gamma)} & = & \pm \frac{1}{2} \left[
     \sqrt{\frac{b_n^{(\alpha, \beta)}}{b_{n + 1}^{(\alpha, \beta)}}} +
     \sqrt{\frac{b_{n - 1}^{(\alpha + 1, \beta + 1)}}{b_n^{(\alpha + 1, \beta
     + 1)}}} \right] .
   \end{array} \]
Using the promotion and demotion three-term recurrences
(\ref{eq:opoly-demote}-\ref{eq:opoly-promote}) we also have the following
recurrence relation:
\begin{equation}
  \label{eq:Psi-sinrecurrence} \begin{array}{lll}
    \Psi_{n + 1}^{(\gamma)} & = & \left[ \tilde{U}_n^{(\gamma)} i \sin \theta
    - \tilde{V}_n^{(\gamma)} \right] \Psi_n^{(\gamma)} + \left[ \tilde{U}_{-
    n}^{(\gamma)} i \sin \theta - \tilde{V}_{- n}^{(\gamma)} \right] \Psi_{-
    n}^{(\gamma)}\\
    &  & \\
    &  & - \tilde{W}_n^{(\gamma)} \Psi_{n - 1}^{(\gamma)} - \tilde{W}_{-
    n}^{(\gamma)} \Psi_{- (n - 1)}^{(\gamma)},
  \end{array}
\end{equation}
where the recurrence constants are given by
\[ \begin{array}{lll}
     \tilde{U}_{\pm n}^{(\gamma)} & = & \frac{1}{2} \left[ \frac{1}{\eta_{n,
     0}^{(\alpha, \beta)}} \mp \frac{1}{\varepsilon_{n - 1, 2}^{(\alpha + 1,
     \beta + 1)}} \right],\\
     &  & \\
     \tilde{V}_{\pm n}^{(\gamma)} & = & \frac{1}{2} \left[
     \frac{\varepsilon_{n - 1, 1}^{(\alpha + 1, \beta + 1)}}{\varepsilon_{n -
     1, 2}^{(\alpha + 1, \beta + 1)}} \pm \frac{\eta_{n, - 1}^{(\alpha,
     \beta)}}{\eta_{n, 0}^{(\alpha, \beta)}} \right],\\
     &  & \\
     \tilde{W}_{\pm n}^{(\gamma)} & = & \frac{1}{2} \left[
     \frac{\varepsilon_{n - 1, 0}^{(\alpha + 1, \beta + 1)}}{\varepsilon_{n -
     1, 2}^{(\alpha + 1, \beta + 1)}} \pm \frac{\eta_{n, - 2}^{(\alpha,
     \beta)}}{\eta_{n, 0}^{(\alpha, \beta)}} \right] .
   \end{array} \]
Finally, putting these last two recurrences together yields
\[ \begin{array}{lll}
     D_n^{(\gamma)} \Psi_{n + 1}^{(\gamma)} & = & \left[ A_n^{(\gamma)} e^{i
     \theta} - B_n^{(\gamma)} \right] \Psi_n^{(\gamma)} + \left[ A_{-
     n}^{(\gamma)} e^{- i \theta} - B_{- n}^{(\gamma)} \right] \Psi_{-
     n}^{(\gamma)} + C_n^{(\gamma)} \Psi_{n - 1}^{(\gamma)} + C_{-
     n}^{(\gamma)} \Psi_{- (n - 1)}^{(\gamma)},
   \end{array} \]
with the following values for the recurrence coefficients:
\[ \begin{array}{lll}
     D_n^{(\gamma)} & = & \left\{ \begin{array}{lll}
       4 \varepsilon_{0, 2}^{(\alpha, \beta)} \sqrt{\gamma}, &  & n = 0,\\
       &  & \\
       2 \varepsilon_{n, 2}^{(\alpha, \beta)} \left[ \sqrt{n + \gamma} +
       \sqrt{n} \right], &  & n > 0.
     \end{array} \right.\\
     &  & \\
     A_{\pm n}^{(\gamma)} & = & \left\{ \begin{array}{lll}
       \sqrt{2} \left[ \sqrt{\gamma + 1} \pm 1 \right], &  & n = 0,\\
       &  & \\
       \sqrt{n + \gamma + 1} \pm \sqrt{n + 1}, &  & n > 0.
     \end{array} \right.\\
     &  & \\
     B_{\pm n}^{(\gamma)} & = & \left\{ \begin{array}{ll}
       \frac{2 \gamma \sqrt{2}}{\sqrt{\gamma + 1}}, & n = 0,\\
       & \\
       \frac{- \varepsilon_{n, 1}^{(\alpha, \beta)}}{2 \sqrt{(n + 1) (n +
       \gamma - 1)}}  \left( \gamma A_{\pm n}^{(\gamma)} + \left[ 2 \sqrt{n (n
       + \gamma)} - 1 \right] A_{\mp n}^{(\gamma)} \right), & n > 0.
     \end{array} \right.\\
     &  & \\
     C_{\pm n}^{(\gamma)} & = & \left\{ \begin{array}{ll}
       0, & n = 0,\\
       & \\
       \frac{1}{A_0^{(\gamma)} (\gamma + 1)} \sqrt{\frac{(2 \gamma + 1)
       \gamma}{\gamma + 2}}, & n = 1,\\
       & \\
       \frac{\gamma \varepsilon_{n, 0}^{(\alpha, \beta)}}{\sqrt{(n + \gamma -
       2) (n + \gamma - 1)} A_{n - 1}^{(\gamma)}} \left[ \sqrt{(n + \gamma)^2
       - 1} - \sqrt{n^2 - 1} \right], & n > 1.
     \end{array} \right.
   \end{array} \]

\section{The Stiffness Matrix}

\label{app:stiff}We assume the decay parameter $s > \frac{1}{2}$ is given and
we derive $\alpha$ and $\beta$ from the value $\gamma \assign s-1$ as in Appendix
\ref{app:recurrence}. Also, we define increments and decrements of the integer
index $k \in \mathbbm{Z}$:
\[ \begin{array}{ccc}
     \alpha \assign - \frac{1}{2}, &  & \beta \assign s - \frac{3}{2},\\
     &  & \\
     k^{\vee} = \tmop{sgn} (k) \left( |k| - 1 \right), &  & k^{\wedge} =
     \tmop{sgn} (k) \left( |k| + 1 \right),\\
     &  & \\
     & n \assign |k| - 1. & 
   \end{array} \]

\noindent We begin by noting the sparse representation of the product of $\phi_k^{(s)}$
and $\frac{1}{(x - i)}$:

\begin{lemma}
  \label{lemma:stiff1}We have the representation:
  \[ \begin{array}{lll}
       \frac{- s}{(x - i)} \phi^{(s)}_k & = & \sum_{l \in \left\{ \pm
       k^{\vee}, \pm k, \pm k^{\wedge} \right\}} \chi_{k, l}^{(s)}
       \phi_l^{(s)},
     \end{array} \]
  for some constants $\chi_{k, l}^{(s)}$.
\end{lemma}

\begin{proof}
  We first note that
  \[ \frac{- s}{x - i} = - \frac{s}{2} \left[ \sin \theta (x) + i (1 + \cos
     \theta (x)) \right], \]
  after making the transformation to $\theta (x)$. Then making the
  identification $\Phi_k^{(s)} = \Psi_k^{(s - 1)}$, we may use recurrence
  relations (\ref{eq:Psi-cosrecurrence})-(\ref{eq:Psi-sinrecurrence}) to obtain
  the result.
\end{proof}

\noindent A second more potent result is the sparsity result for the unweighted Wiener
rational functions $\Phi_k^{(s)} (x)$:

\begin{lemma}
  \label{lemma:stiff2}We have
  \[ \begin{array}{lll}
       \frac{\mathd \Phi_k^{(s)} (x)}{\mathd x} & = & \sum_{l \in \left\{ \pm
       k^{\vee}, \pm k, \pm k^{\wedge} \right\}} \sigma_{k,l}^{(s)} \Phi_l^{(s)},
     \end{array} \]
  with
  \[ \begin{array}{lll}
       \sigma^{(s)}_{k, \pm k^{\vee}} & = & i \tmop{sgn} (k) \frac{n + s}{2 (2
       n + s)} \sqrt{\frac{(2 n + 1) (2 n + 2 s - 1)}{(2 n + s - 1) (2 n + s +
       1)}} \left[ \sqrt{(n + s - 1) (n + s)} \pm \sqrt{n (n + 1)} \right],\\
       &  & \\
       &  & \\
       \sigma_{k, \pm k}^{(s)} & = & i \tmop{sgn} (k) \left\{
       \begin{array}{lll}
         \sqrt{(n + 1) (n + s)}, &  & + k,\\
         &  & \\
         \sqrt{(n + 1) (n + s)}  \frac{s (1 - s)}{(2 n + s) (2 n + s + 2)}, & 
         & - k.
       \end{array} \right.\\
       &  & \\
       \sigma_{k, \pm k^{\wedge}}^{(s)} & = & \pm \frac{n + 1}{n + s + 1}
       \sigma_{k^{\wedge}, \pm k^{}}.
     \end{array} \]
\end{lemma}

\begin{proof}
  This result can be proven by brute-force calculation of the derivatives
  using (\ref{eq:jdemotiona})-(\ref{eq:jdiff}) and the recurrence
  formula (\ref{eq:opoly-3term}). Two critical steps are necessary: a highly
  nontrivial collapsing of a special arithmetic combination involving various
  constants in several Jacobi polynomial relations, and the very special form
  of the $\theta \rightarrow x$ Jacobian for the mapping. Thus, the particular
  form of the mapping is critical in proving this result.
\end{proof}

\noindent Putting the two lemmas together, we have the desired sparsity result for the
$\phi_k^{(s)} (x)$ stiffness matrix:

\begin{theorem}
  \label{thm:stiff-entries}The following equality holds for any $s >
  \frac{1}{2}$:
  \[ \frac{\mathd \phi_k^{(s)}}{\mathd x} = \sum_{l \in \left\{ \pm k^{\vee},
     \pm k, \pm k^{\wedge} \right\}} \tau_{k, l}^{(s)} \phi_l^{(s)}, \]
  where the constants $\tau_{k, l}^{(s)}$ are equal to
  \[ \begin{array}{lll}
       \tau_{k, \pm k^{\vee}} & = & \frac{i}{4} \sqrt{1 - \frac{s (s - 2)}{(2
       n + s - 1) (2 n + s + 1)}} \times\\
       &  & \\
       &  & \left\{ \tmop{sgn} (k) \left( \sqrt{(n + s - 1) (n + s)} \pm
       \sqrt{n (n + 1)} \right) + \right.\\
       &  & \\
       &  & \left. \frac{- s}{(2 n + s)} \left( \sqrt{(n + 1) (n + s - 1)}
       \pm \sqrt{n (n + s)} \right) \right\},\\
       &  & \\
       \tau_{k, k} & = & i \tmop{sgn} (k) \sqrt{(n + 1) (n + s)} - \frac{i
       \nonesep s (s - 1)^2}{2 (2 n + s) (2 n + s + 2)} - \frac{i \nonesep
       s}{2},\\
       &  & \\
       \tau_{k, - k} & = & \frac{i \nonesep s (s - 1)}{2 (2 n + s) (2 n + s +
       2)},\\
       &  & \\
       \tau_{k, \pm k^{\wedge}} & = & \frac{i}{4} \sqrt{1 - \frac{s (s -
       2)}{(2 n + s + 1) (2 n + s + 3)}} \times\\
       &  & \\
       &  & \left\{ - \frac{s}{2 n + s + 2} \left[ \sqrt{(n + 2) (n + s)} \pm
       \sqrt{(n + 1) (n + s + 1)} \right] + \right.\\
       &  & \\
       &  & \left. \tmop{sgn} (k) \left[ \sqrt{(n + 1) (n + 2)} \pm \sqrt{(n
       + s) (n + s + 1)} \right] \right\}.
     \end{array} \]
  Clearly for $|k| = 1$ we have $\tau_{k, \pm k^{\vee}} = \tau_{k, 0}$ so
  that, taking into account the different normalization constant in the
  definition of $\phi_k^{(s)}$ for $|k| = 1$ we have:
  \[ \tau_{k, 0} = \frac{i}{2} \sqrt{\frac{2 s - 1}{2 s + 2}} \left\{
     \tmop{sgn} (k) \sqrt{s} - 1 \right\}, \hspace{1cm} |k| = 1. \]

  \noindent And for $k = 0$:
  \[ \frac{\mathd \phi_0^{(s)}}{\mathd x} = - \frac{i \nonesep \sqrt{s -
     \frac{1}{2}}}{2} \left[ \frac{1 + \sqrt{s}}{\sqrt{1 + s}} \phi_{-
     1}^{(s)} + 2 \sqrt{s - \frac{1}{2}} \phi_0^{(s)} + \frac{1 -
     \sqrt{s}}{\sqrt{1 + s}} \phi_1^{(s)} \right]. \]
\end{theorem}

\begin{proof}
  The hard work was completed in Lemma \ref{lemma:stiff2}. The result of this
  theorem follows from a simple application of the product rule of
  differentiation to
  \[ \frac{\mathd}{\mathd x} \left[ \frac{2^{s / 2}}{(x - i)^s} \Phi_k^{(s)}
     \right], \]
  with the appropriate use of Lemmas \ref{lemma:stiff1} and
  \ref{lemma:stiff2}.
\end{proof}

\noindent With explicit entries for the stiffness matrix derived, we are now ready to
give the proof the third property of Theorem \ref{thm:stiffness-brief}, the
spectral radius of the stiffness matrix.

\begin{proof}
  We can crudely bound the entries of the stiffness matrix for $n \assign |k|
  - 1 > 1$
  \[ \begin{array}{lll}
       | \tau_{k, k^{\vee}} | & \leq & \frac{n}{2} + s,\\
       &  & \\
       | \tau_{k, - k^{\vee}} | & \leq & \frac{s}{2},\\
       &  & \\
       | \tau_{k, \pm k^{\vee}} | & \leq & \frac{n}{2} + s,\\
       &  & \\
       | \tau_{k, k} | & \leq & n + 2 s,\\
       &  & \\
       | \tau_{k, - k} | & \leq & 1,\\
       &  & \\
       | \tau_{k, k^{\wedge}} | & \leq & \frac{n}{2} + s + \frac{1}{2},\\
       &  & \\
       | \tau_{k, - k^{\wedge}} | & \leq & \frac{s}{4}.
     \end{array} \]

  \noindent Thus we have:
  \[ \begin{array}{lll}
       | \tau_{- k} | + | \tau_{k^{\vee}} | + | \tau_{- k^{\vee}} | + |
       \tau_{k^{\wedge}} | + | \tau_{- k^{\wedge}} | & \leq & 1 + \frac{n}{2}
       + s + \frac{s}{2} + \frac{n}{2} + s + \frac{1}{2} + \frac{s}{4}\\
       &  & \\
       & = & n + 3 s + 2.
     \end{array} \]
  An application of Gerschgorin's Theorem then proves the result.
\end{proof}

%% file: wienerbasis-part1.bbl
\begin{thebibliography}{10}
\bibitem{whittaker1915}
E.T. Whittaker, \emph{On the functions which are represented by the expansions
  of the interpolation theory}, Proceedings of the Royal Society of Edinburgh
  \textbf{35} (1915), 181--194.
\end{thebibliography}
